\documentclass[12pt]{amsart}
\usepackage{amsfonts}
\usepackage{amsmath}
\usepackage{amsxtra}
\usepackage{amssymb,latexsym}
\usepackage[mathcal]{eucal}
\usepackage{amscd}

\newtheorem{theo}{{\bfseries Theorem}}[section]
\newtheorem{prop}[theo]{{\bfseries Proposition}}
\newtheorem{lem}[theo]{{\bfseries Lemma}}
\newtheorem{cor}[theo]{{\bfseries Corollary}}
\newtheorem{df}[theo]{{\bfseries Definition}}

\def \C {\mathbb C}

\def \Z {\mathbb Z}
\def \R {\mathbb R}

\def \E {\mathcal E}

\def \ep {\epsilon}

\def \om {\omega}
\def \r {\rho}

\def \ii {\mathbf{i}}

\usepackage{amsfonts}
\usepackage{amsmath}
\usepackage{amsxtra}
\usepackage{amssymb,latexsym}
\usepackage[mathcal]{eucal}
\usepackage{amscd}

\numberwithin{equation}{section}

\begin{document}
\begin{titlepage}

\title{The Simple Yield Curve Models}

\author{Ethan Akin and Morton Davis}
\address{Mathematics Department,
 The City College, 137 Street and Convent Avenue,
 New York City, NY 10031, USA}

\email{ethanakin@earthlink.net, mortdavis@att.net}

\date{September 2015, March 2024}

\begin{abstract} With $P_t$ the price in current dollars of  a dollar delivered $t$ time units from now, we assume that
 $P$ is a decreasing function defined for $t \in \R_+$ with $P_0 = 1$. The  negative logarithmic derivative, $- \stackrel{\bullet}{P}_t/P_t$
 defines the yield curve function $Y_t$. An $n$ parameter linear yield curve model selects as allowable yield curves
 $Y_t(r) = \sum_{i=1}^n r_i Y^i_t$ with the functions $Y^i$ fixed and with $r$ varying over an open subset of $\R^n$ on which
 $Y_t(r) \ge 0$ for all $t \in \R_+$. For example, the flat yield curve model with $P_t(r) = e^{-rt}$ is a one parameter linear model with
 $Y^1_t(r) = r > 0$. We impose two natural economic requirements on the model: (SPA) static prices allowed, i.e. it is always possible that
 as time moves forward, relative prices do not change, and (NLA) no local arbitrage, i.e. there does not exist a self-financing bundle of futures
 such that the zero present value is a local minimum with respect to small changes in the space of admissible yield curves.  In that case
 the model always contains one of four simple models.  If we impose the additional requirement (LRE) long rates exist, i.e. for every $r$
  $Lim_{t \to \infty} Y_t(r)$ exists as a finite limit, then the number of simple models is reduced to two.
\end{abstract}



\end{titlepage}
 \maketitle

\tableofcontents
\setcounter{page}{1}

\section{Introduction}

Acknowledgement: We would like to express our thanks to our colleague Eli Amzallag for helpful discussions about this work.
 \vspace{.5cm}

 Notation: As usual, we use $\R$ and $\C$ for the real and complex fields and $\Z$ for the ring of integers. We use $\R_+$ to denote the
 set of non-negative reals. That is, $\R_+ = [0,\infty)$.
  \vspace{1cm}

When the interest rate is $r$, dollars delivered $t$ years in the
future are discounted by the factor $e^{-rt}$.  To put this
another way the price $P_t$ of one dollar delivered $t$ years
from now is $e^{-rt}$ in current dollars.  Here $t$ varies
continuously in  $\R_+$.  Thus, the yield curve
$Y_t$, defined to be the negative logarithmic derivative of
$P_t$, i.e. $Y_t = - \stackrel{\bullet}{P}_t/P_t$, is the
constant $r$.  Even if we ignore risk, which induces a price
spread among $t$ futures for fixed $t$, the price curve $P_t$ may
not exhibit this simple shape.  To describe the pattern of
interest rates we then need the entire yield curve $Y_t$ instead
of the single number $r$.

Nonetheless, when we think of interest rates as unchanging, the
flat yield curve model with its single interest rate is a simple
and useful tool adequate for many theoretical purposes.  However,
when we begin to consider change in interest rates, the flat
model won't do.  This paper had its origin in our discovery of
its flaws.

We had, in fact, rediscovered an odd phenomenon well-known among
financial economists who consider the notion of {\em duration}
for a bond.

Consider a finite set $D$ consisting of several future times.
For each time $t \in D$ allot $z_t > 0$ current dollars for $t$
futures ($z_t$ is defined to be 0 for $t \not\in D$).  If the
current interest rate is $r^*$, this allocation buys a bond
paying $x_t$ dollars $t$ years from now where
\begin{equation}\label{1.01}
z_t \ = \ e^{-r^*t}x_t.
\end{equation}

The present value of the basket of futures is
\begin{equation}\label{1.02}
V_0 \ = \sum_t \ z_t \ = \sum_t \ e^{-r^*t}x_t.
\end{equation}
Following \cite{M} (see also \cite{H}) we
define the duration $T$ of the basket of futures as the weighted
average of the payment times:
\begin{equation}\label{1.03}
\begin{split}
T \ = (\sum_t \ tz_t)/(\sum_t \ z_t) \hspace{1cm}\\
= \ (\sum_t \ t  e^{-
tr^*}x_t)/(\sum_t \ e^{-tr^*}x_t).
\end{split}
\end{equation}
Alternatively, we can define the duration as the unique
time $T$ satisfying the -equivalent- equations
\begin{align}\begin{split}\label{1.04}
0 \ = &\sum_t \ (T - t)z_t, \\
0 \ = &\sum_t \ (T - t)e^{(T-t)r^*}x_t.
\end{split}\end{align}

Suppose we finance this bond by a single sale of $V_0$ worth of
time $T$ futures.  That is, define
\begin{equation}\label{1.05}
\tilde{x}_s \ = \left \{ \begin{array}{cr}
x_s & s \neq T\\
x_T \ - \ \sum_t \ e^{(T-t)r^*}x_t & s = T
\end{array} \right.
\end{equation}
or, equivalently,
\begin{equation}\label{1.08}
\tilde{z}_s \ = \left\{ \begin{array} {lr}
z_s & s \neq T\\
z_T - (\sum_tz_t) & s = T,
\end{array} \right.
\end{equation}

The investment $\tilde{x}_t$ is self-financing, i.e.
its present value, $\tilde{V}_0$ is zero.

Fisher and Weil in \cite{FW} point out that this sale immunizes the
investor against loss when interest rates change.  To see this,
suppose that a period $h$ as passed, short enough that none of
the payments are due at or before $h$, and that the new interest
rate is now $r$.  The quantity $\tilde{x}_t$ is a payment
occurring $t - h$ in the future.  So the new present value is:
\begin{equation}\label{1.06}
\begin{split}
\tilde{V}_h \ = \sum_t \ e^{-(t-h)r} \tilde{x}_t \ = \hspace{2cm} \\
e^{-(T-h)r}[\sum_t \ e^{(T-t)r}x_t - \sum_t \ e^{(T-t)r^*} x_t].
\end{split}
\end{equation}

By applying (\ref{1.04}) we can write this as
\begin{equation}\label{1.07}
\tilde{V}_h \ = e^{-(T-h)r} [\sum_t \ e^{(T -
t)r^*}x_t[e^{(T-t)(r-r^*)}\ - \ 1 \ - \ (T-t)(r-r^*)].
\end{equation}

Observe that $e^a - 1 - a = \int^a_0(e^s - 1)ds$ is
nonnegative for all $a$, vanishing only when $a = 0$.  Thus, we
obtain $\tilde{V}_h \geq 0$ vanishing only when $r = r^*$.

As Ingersoll, Shelton and Weil remark, \cite{ISW} page 635, this
result is ``too good''.  It is an arbitrage opportunity which
exists at any interest rate.  In fact, (\ref{1.03}) shows that $T$ can
be defined using the current dollar quantities $z_t$ in a way
independent of the interest rate.  

Arbitrage opportunities do happen.  Noticing such an investment,
a speculator can reap the associated profit.  However, the
transactions involved cause compensating price movements which
tend to  wipe out the original opportunity.  Here, however, as
long as the yield curve remains flat and will continue to do so
the arbitrage possibilities described above remain as well.  Such
speculations in the futures market would, instead, bend the yield
curve.

A yield curve model selects as admissible certain special curves
together with their associated price curves.  For example, for
the flat yield curve model the positive constants are the
admissible yield curves.  An arbitrage opportunity for a model is
an investment with no possibility of loss when the prices curves
are restricted to those admissible according to the model.  Our
first condition on a yield curve model is that even local
versions of such arbitrages do not occur.

The second condition is that the model allows static prices.
That is, starting from an admissible price curve the possibility
that relative prices among futures do not change is not excluded
from the menu of choices offered by the family of admissible
price curves.  In the next section, we will examine in some
detail what this means.

Our third condition is linearity.  This is a mathematical
convenience rather than an economic imperative.  The flat yield
curve model is described by $Y_t(r) = r \cdot 1$ with a single
positive parameter.  An $n$ parameter linear model is of the form
$Y_t(r_1, \dots, r_n) = r_1 Y^1_t + \dots + r_n Y^n_t$ where
$(r_1, \dots, r_n)$ varies over some open subset of $\R^n$
and the curves $Y^1_t, \dots, Y^n_t$ are fixed for the model.

The infinitesimal or short rate for the yield curve $Y_t$ is
$Y_0$.  The limit or long rate is $Lim_{t \rightarrow \infty}\ Y_t$
if the limit exists.  We say the model admits long rates if this
limit exists as a finite real number for any admissible yield
curve.

Our main result is that an $n$ parameter linear model allows
static prices, admits no local arbitrage, contains the flat yield
curves and admits long rates if and only if it contains, in a
natural way to be described later, one of the following two
models:
\begin{equation}\label{1.09}
Y_t(r_1,r_2,r_3) \ = r_1 + r_2 e^{-\r t}
+ r_3 e^{-2 \r t}
\end{equation}
where $\r$ is a fixed positive constant, or
\begin{equation}\label{1.10}
Y_t(r_1,r_2,r_3,r_4) \ = r_1 + r_2 e^{-
\r t}\cos \om t + r_3 e^{-\r t}\sin \om t + r_4 e^{-2\r t}
\end{equation}
where $\r$ and $\omega$ are fixed positive constants.

In each case $r_1$ is the long rate.  If we further demand that
the approach to the long rate be eventually monotonic then the
damped oscillation of (\ref{1.10}) is excluded and we are left with
(\ref{1.09}).

If we omit the long rate demand and allow the yield curves to be
unbounded then the constant $\r$ in (\ref{1.09}) can be negative, and
we have two further possible models:
\newpage
%
\begin{equation}\label{1.11}
Y_t(r_1,r_2) \ = r_1 + r_2 t, \hspace{6cm}\end{equation}
$$ Y_t(r_1,r_2,r_3,r_4,r_5) \ = r_1 + r_2 \cos
\om_1 t + r_3 \sin \om_1 t + r_4 \cos \om_2 t + r_4 \sin \om_2 t.$$
 where $\om_1, \om_2$ are positive constants such
that $\om_1/\om_2$ is irrational.

A wide variety of yield curve models arise as diffusion models.
Each of these is an $n$ parameter linear model for some $n$.
There is an additional dynamic structure in that the change in
prices over time is governed by a diffusion stochastic process in
the parameter space $\R^n$.

This is not how diffusion models are usually presented.  Instead
(e.g. see Ingersoll, Shelton and Weil (1978) pages 637ff.) from
diffusion equations for general assets the yield curve at each
parameter point is obtained by solving an associated partial
differential equation.  The result is a linear yield curve model
with no local arbitrage.

These models do not usually satisfy the static price condition.
This is not serious because the static price condition can always
be obtained by introducing an additional translation parameter
(see Section 1 below).  However, in most diffusion models the new
parameter enters in a nonlinear fashion.  Thus, the new enlarged
model no longer satisfies linearity.

It is the interplay between these three conditions which limits
the range of possible models.  If we impose the additional
condition that each yield curve approaches a nonzero long rate
limit in a monotone fashion then by rescaling time we can assume
$\r$ in (\ref{1.09}) is $1$.  We thus obtain the unique simple yield
curve model, the three parameter linear model
\begin{equation}\label{1.12}
Y_t(r_1,r_2,r_3) \ = r_1 + r_2 e^{-t} + r_3 e^{-2t},
\end{equation}
where $(r_1,r_2,r_3)$ varies over the open set in $\R^3$ such that
\begin{equation}\label{1.13}
inf_{t \geq 0} \ \
Y_t(r_1,r_2,r_3) \quad > \quad 0.
\end{equation}
\vspace{1cm}

\section{Yield Curve Models}

In general, \emph{price} describes the exchange ratio between
quantities of two commodities in a market.  If among many
commodities we choose one, say corn, as \emph{numeraire} then the
price $p_i$ of commodity $i$ is an amount of corn for which a
unit quantity of $i$ exchanges.  A basket of commodities is
described by a vector $x$ such that coordinate $x_i$ is the
quantity of $i$ in the basket.  The cost or \emph{value} of such a
basket, given by
\begin{equation}\label{2.01}
V(x) = \sum_i \ p_i x_i,
\end{equation}
is the quantity of corn for which the entire basket can be
exchanged.  In (\ref{2.01}) $x_i$ can be positive or negative.  We
interpret a negative $x_i$ as a liability, a commitment to
deliver quantity $|x_i|$ of good $i$.  So $V(x)$ is the net value
summed over the goods coming in and going out.  To change
numeraire from corn to gold, observe that if $p_g$ is the corn
price of gold then $p_i/p_g$ is the gold price of commodity $i$.
So we get the gold value of basket $x$ by dividing $V(x)$ by
$p_g$.

We assume that all prices are nonnegative.  This is a consequence
of the (often tacit) assumption of \emph{free disposal}.
Similarly, the assumption that choice, as well as disposal, is
costless implies that the market evaluates a choice between goods
$A$ and $B$ as at least as valuable as good $A$ alone.  For
mathematical convenience we will assume that prices are nonzero
as well and so are strictly positive.

We are also ignoring transaction costs and are assuming that there is a single price
to buy or to sell a commodity in the market.

The description of a commodity requires a specification of time
and place of delivery.  We will ignore place but for futures
markets the date of delivery is exactly the central concern.
Thus, the \emph{price curve} for dollar futures is defined by:
\begin{quote}
 $P_t = $ the current dollar price of one \\
dollar delivered $ t $  time units from now.
\end{quote}

The \emph{present value} of a finite bundle of payments $x = \{x_t\}$ is
the value in current dollars computed as in (\ref{2.01}).
\begin{equation}\label{2.02}
V(x) \ = \sum_t \ P_t x_t.
\end{equation}
As before $x_t$ is positive or negative according to whether the
quantity $|x_t|$ of dollars is incoming or outgoing $t$ time
units from now.

The present value is a sum rather than an integral because $x$
consists of a finite number of discrete payments rather than a
continuous flow.  For any set $T$, we call $x: T \rightarrow \R$ a \emph{finite function}
when the \emph{support} of $x$,
denoted supp$(x)$, is a finite subset of $T$, where
\begin{equation}\label{2.03}
{\rm supp}(x) \ = \{ \ t \ : \ x_t \neq 0 \ \}.
\end{equation}
So a bundle of payments $x$ is a finite function on $\R_+$.
We call $x$ a {\em bundle of futures} when $x_0 = 0$, i.e. $0
\not\in {\rm supp}(x)$, so that no current payments are involved.

Observe that we use continuous variables for time and later for
interest rate parameters.  In contrast with our finite payment
functions, the price functions are assumed to be {\em smooth}.
That is, continuous partial derivatives of all orders exist.

The numeraire for our price function $P$ is current dollars and
so $P_0 = 1$, the current dollar price of one current dollar.
Furthermore, the price function is nonincreasing.  This follows
from the free disposal assumption and the further assumption of
\emph{free storage.}  To each consumer a dollar delivered at an
early date is at least as valuable as a later one because one
could always store the former and so obtain the latter at no
additional cost.  In other words, if we regard the late dollar as
good $A$ then the early dollar represents a choice between good
$A$, obtained by storage, and good $B$, the consumption of the
dollar during the intervening time.  In particular, writing
$\stackrel{\bullet}{P}_t$ for the derivative at time $t$, we have
that $\stackrel{\bullet}{P}_t$ is nonpositive.  We will assume it
is always negative.  Thus, we have the following conditions on
the price curves, defined for all $t \geq 0$:
\begin{equation}\label{2.04}
\begin{split}
1 \geq P_t > 0 \qquad \mbox{ with} \quad P_0 = 1, \\
\stackrel{\bullet}{P}_t < 0. \hspace{3cm}
\end{split}
\end{equation}

Associated with the price function there are several related
curves.  The \emph{log-price} is defined by
\begin{equation}\label{2.05}
L_t \ = \ - \ln P_t,
\end{equation}
(where ln is the natural log) and its derivative, the \emph{yield
curve} is defined by
\begin{equation}\label{2.06}
Y_t \ = \ \stackrel{\bullet}{L}_t \ = \ - \stackrel{\bullet}{P}_t/P_t.
\end{equation}
These satisfy the conditions
\begin{align}\label{2.07}
\begin{split}
&L_t \geq 0 \quad \mbox{with} \quad L_0 \ = \ 0, \\
&Y_t > 0.
\end{split}
\end{align}

From the initial conditions $P_0 = 1$, or equivalently $L_0 =
0$, we can recover the price from the yield curve
\begin{align}\label{2.08}
\begin{split}
&L_t \ = \int^t_0 \ Y_s ds,\\
&P_t \ = \exp(-L_t).
\end{split}
\end{align}

Observe that we are using the phrase ``yield curve'' for the
marginal value of the log-price.  It is sometimes used instead to
refer to the average value $L_t/t$.  The two concepts agree
exactly when the yield curve is a constant.  In this, the \emph{flat yield curve} case,
with $Y_t = r > 0$, we have $L_t = rt$
and $P_t = \exp(-rt)$.  This is the classical price curve when
``the interest rate'' is $r$.

A \emph{yield curve model} picks out from the set of all smooth, positive real
valued functions on $\R_+ = [0,\infty)$ a subset, the \emph{admissible yield curves},
with associated price curves given by (\ref{2.08}).  An
$n$ \emph{parameter yield curve model} associates to each
admissible yield curve a unique vector $r$ in some subset $O$ of
$\R^n$.  To be more precise, such a model is a function $Y$
from $O$ to the vector space $\R^{\R_+}$ of real
valued functions on $\R_+$, with special assumptions as
follows:
\vspace{.5cm}

\begin{df}\label{df2.01}  Let $n$ be a positive integer
and $O$ be a nonempty open subset of $\R^n$.  An \emph{ $n$
parameter yield curve model defined on $O$} is a smooth, positive
real valued function in $n + 1$ real variables, $Y_t(r_1, \ldots,
r_n)$ defined for all $t \geq 0$ and $r = (r_1, \dots , r_n)$ in
$O$.  Regarded as a function from $O$ to $\R^{\R_+}$,
we assume that $Y$ is injective.  That is, if for $r,r' \in
O$, $Y_t(r) = Y_t(r')$ for all $t \in \R_+$, then $r = r'$.  \end{df}
\vspace{.5cm}

For example, the \emph{flat yield curve model} is the one
parameter model $Y_t(r_1) = r_1$ defined on $O = (0,\infty)$.

Changing prices is represented as continuous motion in the set of
admissible yield curves.  For an $n$ parameter model defined on
$O$, such change in time is described by smooth motion of the
parameter $r$ in the set $O$.

We now describe the conditions which we will require of a good
model.
\vspace{.5cm}

\subsection{Static Prices Allowed (SPA) }

If $Y_t$ is an admissible
yield curve whose associated price curve $P_t$ represents prices
in the current market, then we want the model to include the
possibility that over time relative prices do not change.  We are
not claiming that this sort of stasis is at all likely, but we do
not want it ruled out a priori by the model.
\vspace{.5cm}

At first glance this condition appears to be trivially satisfied.
We are fixing the set of admissible curves in a time independent
way so if $Y_t$ is admissible today, then it is admissible
tomorrow.  However, an unchanging yield curve is usually not the
same thing as unchanging prices.  The confusion between these two
concepts comes from the two different roles played by time in the
futures market.

On the one hand, each instant on the time line is a \emph{date},
$d$, a moment on the calendar.  The date is the label attached to
our commodities.  For precision we will use the notation
\begin{equation}\label{2.09}
\$_d \ \equiv \  \mbox{ a dollar  delivered  at  date } \ d.
\end{equation}

On the other hand, between two dates is an interval of time.
Suppose that from a date $d_0$ we move to a date $h$ time units
later which we will denote $d_h$.  Then,
\begin{equation}\label{2.10}
d_h \ = d_0 + h \quad \mbox{and} \quad h \ = d_h - d_0.
\end{equation}
The variable $t$ in the price and yield curves is this interval
length measuring the wait from now, $d_0$, to a later date $d_t = d_0 + t$.
Thus, if $P$ is the price function now at $d_0$, then $P_t$ is
the price of $\$_{d_0+t}$ in current dollars, i.e. in $\$_{d_0}$.

Suppose that when we move from the date $d_0$ to the date $d_h$
the price, log-price and yield functions change from $P$, $L$,
$Y$ to new ones $P^h$, $L^h$, $Y^h$.

The numeraire has changed from $\$_{d_0}$ to $\$_{d_h}$.  $P^h_t$
is the price of $\$_{d_h+t}$ in the new numeraire $\$_{d_h}$.  In
the old price system the commodity $\$_{d_h + t} = \$_{d_0 + h
+t}$ is priced at $P_{t+h}$ valued in $d_0$ dollars, $\$_{d_0}$.

Similarly, suppose that we purchased at date $d_0$ a bundle of
futures $x$ and that no payments have occurred in the interval
from $d_0$ to $d_h$, that is
\begin{equation}\label{2.10a}
x_t \ = 0 \quad \mbox{for} \quad 0 \leq t \leq h
\end{equation}
or equivalently, $[0,h] \cap supp(x) \ = \emptyset.$

For each $t$, $x$ describes the quantity $x_t$ of commodity
$\$_{d_0 + t}$.  At date $d_h$ the asset remains the same but is
described by the new finite function $x^h$ with
\begin{equation}\label{2.11}
x^h_{t-h} \ = x_t,
\end{equation}
that is, the wait for each commodity has been reduced by $h$.  The
$d_h$ present value is now
\begin{equation}\label{2.12}
V^h \ = \sum_t \ P^h_t x^h_t \ = \sum_t \ P^h_{t-h}x_t.
\end{equation}

Notice that $x_t = 0$ for $0 \le t < h$ implies that although $P^h_{t-h}$ is undefined for
$t$ in this range we can still regard both sums with $t$ varying in $\R_+$.

\emph{Unchanging prices} means that the relative prices of goods $\$_d$
with $d \geq d_h$ are the same with respect to the price function
$P^h$ as they were with respect to $P$.  That is, for all $t \geq 0$
\begin{equation}\label{2.13}
\begin{split}
P^h_t \ = P_{t+h}/P_h, \hspace{3cm}\\
L^h_t \ = L_{t+h} - L_h, \hspace{3cm} \\
Y^h_t \ = Y_{t+h}. \hspace{3.5cm}
\end{split}
\end{equation}

The first equation describes unchanging relative prices.  We are merely changing the numeraire from $\$_{d_0}$ to $\$_{d_h}$. The
second equation follows by taking logs and the third, in turn, by taking the derivative with respect to $t$ because the
derivative of the constant $-L_h$ is zero.  On the other hand, by
using (\ref{2.08}) we can recover the first, price, equation from the
third, yield, equation in (\ref{2.13}).

Thus, the SPA condition says that if $Y$ is an admissible curve,
then for any $h > 0$ the $h$ translated curve $\Phi_h(Y)$ is also
admissible where
\begin{equation}\label{2.14}
\Phi_h(Y)_t \ \equiv \ Y_{t+h} \qquad t,h \geq 0.
\end{equation}

For an $n$ parameter model defined on $O$ the $h$ translate of
the curve $Y(r)$ for $r \in O$ must then correspond to a unique
parameter value which we denote $\varphi_h(r)$.  That is,
\begin{equation}\label{2.15}
\Phi_h(Y(r)) \ = Y(\varphi_h(r))
\end{equation}
We assume the function $\varphi$ is smooth.  Thus, we are led to
\vspace{.5cm}

\begin{df}\label{df2.02}  Let $Y_t(r)$ be an $n$
parameter yield curve model defined for $(t,r) \in \R_+ \times O \subset \R^{n+1}$,
with associated log-price and
price curves:  $L_t(r) \ = \int^t_0 Y_s(r)ds$ and $P_t(r) \ = \exp[-L_t(r)]$.
We say that the model {\em allows static
prices}, or satisfies condition SPA, if there is a smooth
function $\varphi: \R_+ \times O \rightarrow O$ such that
for all $(t,h,r) \in \R_+  \times \R_+  \times O$
\begin{equation}\label{2.16}
\begin{split}
P_t(\varphi_h(r)) \ = P_{t+h}(r)/P_h(r) \hspace{2cm}\\
L_t(\varphi_h(r)) \ = L_{t+h}(r) - L_h(r) \hspace{2cm}\\
Y_t(\varphi_h(r)) \ = Y_{t+h}(r).\hspace{2.5cm}
\end{split}
\end{equation}
\end{df}
\vspace{.5cm}

{\bf Remark.}  As with (\ref{2.13}) the three equations in (\ref{2.16}) are
equivalent.  By Definition \ref{df2.01} the association from parameters to
yield curves is assumed injective and so the parameter value
$\varphi_h(r)$ is uniquely determined for each $(h,r) \in
\R_+ \times O$ by the yield curve equation in $t$ of
(\ref{2.16}).

\vspace{.5cm}

From the remark above we obtain the Kolmogorov flow conditions
\begin{equation}\label{2.17}
\begin{split}
\varphi_0(r) \ = r; \hspace{3cm} \\
\varphi_{h_1}(\varphi_{h_2}(r)) \ = \varphi_{h_1 + h_2}(r) \hspace{1cm}
\end{split}
\end{equation}
for all $r \in O$ and $h_1,h_2 \in \R_+ $. Observe that the two
sides of each equation have the same yield curve.  We will call
the function $\varphi$ the {\em static price flow} on $O$.

Notice that if we define the {\em short rate function} $Y_0: O
\rightarrow (0,\infty)$ by setting $t = 0$ in $Y_t(r)$ then we
can recover the full yield curve from the short rate function and
the static price flow $\varphi$ by using (\ref{2.16}):
\begin{equation}\label{2.18}
Y_t(r) \ = Y_0(\varphi_t(r)) \hspace{2cm}
\end{equation}
for all $(t,r) \in \R_+  \times O$.

Translation of a function by $h$ leaves the function unchanged
for all positive $h$ if and only if the function is constant.
Thus, the flat yield curve is the one case where the notions of
unchanging yield curve and unchanging prices agree.  In
particular, the one parameter flat yield curve model allows
static prices with $\varphi_h(r_1) = r_1$ for all $(h,r_1) \in
\R_+  \times (0,\infty)$.

A smooth flow like $\varphi$ is the solution of a smooth system
of ordinary differential equations.  For each $r \in O$ let
$\stackrel{\bullet}{\varphi}_h(r)$ denote the derivative $d \varphi_h(r)/dh$, i.e.
the velocity vector of the path in $O$,
$\varphi_h(r)$, obtained by fixing $r$ and letting $h$ vary.  The
vector field $\xi: O \rightarrow \R^n$ is obtained by taking
the velocity vector at $h = 0$, that is,
\begin{equation}\label{2.19}
\xi(r) \ = \ \stackrel{\bullet}{\varphi}_0(r).
\end{equation}
Differentiating (\ref{2.17}) with respect to $h_1$ and setting $h_2 =
h$, we see that $\varphi_h(r)$ is the solution of the smooth
initial value problem
\begin{equation}\label{2.20}
\begin{split}
\varphi_0(r) \ = r \hspace{2cm} \\
\stackrel{\bullet}{\varphi}_h(r) \ = \ \xi(\varphi_h(r)). \hspace{1cm}
\end{split}
\end{equation}
Definition \ref{df2.02} includes the assumption that this solution path is
defined for all positive $h$.  In o.d.e. terminology this says
that the vector field $\xi$ is \emph{complete} in the positive
direction.

With $O$  an open set, we can use this characterization of the
flow to extend $\varphi_h(r)$ to some negative values of $h$.

\vspace{2ex}

\begin{lem}\label{lem2.03} Let $Y_t(r)$ be an $n$ parameter
yield curve model defined on an open subset $O$ of $\R^n$,
with $P_t(r)$ and $L_t(r)$ the associated price and log-price
curves.  Assume that the model satisfies SPA with $\varphi:
\R_+  \times O \rightarrow O$ the static price flow.  There
exists an open subset $\tilde{O}$ of $\R \times O$ which
contains $\R_+  \times O$, and a smooth function
$\tilde{\varphi}: \tilde{O} \rightarrow O$ which extends
$\varphi$, and for $(t,-h,r) \in \R_+  \times \tilde{O}$
such that $h,t - h \geq 0$
\begin{equation}\label{2.21}
\begin{split}
P_{t-h}(r) \ = P_t(\tilde{\varphi}_{- h}(r))/P_h(\tilde{\varphi}_{-h}(r)) \\
L_{t-h}(r) \ = L_t(\tilde{\varphi}_{-h}(r)) - L_h(\tilde{\varphi}_{-h}(r)) \\
Y_{t-h}(r) \ = Y_t(\tilde{\varphi}_{-h}(r)). \hspace{2cm}
\end{split}
\end{equation}
\end{lem}

\proof  The constructions and proof of the existence and
uniqueness theorem for o.d.e.'s yield the following results, see,
for example, \cite{L} Chapter IV Sections 1 and 2.  The
solution flow for the vectorfield $\xi$ is a smooth function
$\tilde{\varphi}: \tilde{O} \rightarrow O$ with $\tilde{O}$ an
open subset of $\R \times O$.  For each $r \in O$, $\{h:
(h,r) \in \tilde{O}\}$ is an open interval containing $0$.  It is
the maximum open interval on which the unique solution path of
the initial value problem (\ref{2.20}) is defined.  This interval thus
includes $\R_+ $ because $\xi$ is assumed positively
complete.  Furthermore, by uniqueness, $\tilde{\varphi}$ extends
$\varphi$.  Finally, the flow equation (\ref{2.17}) for
$\tilde{\varphi}$ holds in the sense that if the left side is
defined then the right side is as well.

In particular, if $(t,-h,r)$ is in $\R_+  \times \tilde{O}$
with $h,t-h \geq 0$, then $\varphi_h(\tilde{\varphi}_{-h}(r)) =
r$.  Now let $\tilde{r} = \tilde{\varphi}_{-h}(r)$, $\tilde{t} = t - h$
and apply equations (\ref{2.16}) to $(\tilde{t},h,\tilde{r})$.  Since
$\tilde{t} + h = t$ and $\varphi_h(\tilde{r}) = r,$ we get
(\ref{2.21}).

$ \Box$\vspace{.5cm}

Thus, condition SPA imposes some structure upon an $n$ parameter
model.  On the other hand, Joe Langsam pointed out that we can
always obtain SPA by throwing in another parameter.
\begin{equation}\label{2.22}
\begin{split}
\tilde{Y}_t(r_1, \dots, r_n,r_{n+1}) \ = Y_{t+r_{n+1}}(r_1,\dots, r_n) \\
\tilde{\varphi}(r_1, \dots, r_n,r_{n+1},h) \ = (r_1, \dots, r_n,r_{n+1} + h)
\end{split}
\end{equation}
is an $n+1$ parameter model on $O \times \R_+$ with
$\tilde{\varphi}$ the required static price flow.
\vspace{.5cm}

%

\subsection{No Local Arbitrage (NLA)}
\vspace{.5cm}

A bundle of futures $x$ is
called a \emph{self-financing investment} if at current prices the
present value is zero.  This means that the worth of the assets
and liabilities in the bundle exactly balance and no current cash
is required for $x$.  Since we are ignoring transaction costs,
anyone can obtain $x$ in the current market.

Suppose that $h$ is small enough that no payments occur during
the initial $h$ time interval, i.e. condition (\ref{2.10a}) holds.  The
bundle $x$, now labeled $x^h$, is worth $V^h$ computed using the
new price curve $P^h_t$, cf. equations (\ref{2.11}) and (\ref{2.12}).  For
example, if relative prices are unchanged so that $P^h$ is
related to $P$ by equation (\ref{2.13}) then the value $V^h$ is still
zero. We have merely changed the numeraire from $\$_{d_0}$ to $\$_{d_h}$.
In general, when prices have changed $V^h$ is pure profit
when positive, and loss when negative, measured in now current
dollars $\$_{d_h}$.

For a yield curve model a bundle of futures $x$ is an arbitrage
for an admissible price function $P$ if it has present value zero
with respect to $P$ and if for any $h$ satisfying (\ref{2.10a}) and for
any admissible price function $P^h$ the value $V^h$ given by
(\ref{2.12}) is nonnegative.  Thus, $x$ is a self-financing investment
at current prices and no admissible pattern of price movements
returns a loss on the investment.  As described in the
Introduction the occurrence of such arbitrage possibilities is a
serious flaw in a yield curve model.  We want to exclude even the
local versions of such arbitrage.
\vspace{.5cm}

\begin{df}\label{df2.04} Let $Y_t(r)$ be an $n$
parameter yield curve model defined on $O \subset \R^n$ with
associated price curve $P_t(r) \ = \exp[- \int^t_0 Y_s(r)ds]$.  Let
$x$ be a bundle of futures, i.e. a finite function $x: \R_+
\rightarrow \R$ with $x_0 = 0$.  $x$ is called {\em self-financing at} $r \in O$ if
\begin{equation}\label{2.23}
\sum_t \ P_t(r)x_t \ = 0.
\end{equation}
$x$ is called a {\em null investment at} $r$ if for all
$(h,r^\prime)$ in some neighborhood $[0,\ep) \times U$ of $(0,r)$ in $\R_+ \times O$ with $[0,\ep] \cap {\rm supp}(x) = \emptyset$
\begin{equation}\label{2.24}
\sum_t \ P_{t-h}(r^\prime)x_t \ = 0.
\end{equation}
$x$ is called a {\em local arbitrage at} $r$ if $x$ is self-financing
at $r$ and if for all $(h,r^\prime)$ in some
neighborhood $[0,\ep) \times U$ of $(0,r)$ in $\R_+ \times O$ with $[0,\ep] \cap {\rm supp}(x) = \emptyset$
\begin{equation}\label{2.25}
\sum_t \ P_{t-h}(r^\prime) x_t \ \geq \ 0.
\end{equation}
We say that the model allows {\em no local arbitrage}, or
satisfies condition NLA, if the zero function, $x_t = 0$ for all
$t \in \R_+ $, is the only local arbitrage at any $r \in
O$. \end{df}
\vspace{.5cm}

Thus, condition NLA has two aspects.  First, it requires that $0$
is the only null investment.  This means that any two different
bundles of futures can be distinguished by some price movement
near $r$ in the market.  Notice that two bundles of futures have
the same value not only now with respect to $P_t(r)$ but at time
$h$ from now at $P_t(r^\prime)$, for all $(h,r^\prime)$ near
$(0,r)$, exactly when the difference between the two bundles is
null.  Secondly, the condition excludes entirely true arbitrage
opportunities, investments self-financing now for which no local
motion leads to a loss and some local motion leads to a profit.

When the model satisfies condition SPA the description of a local
arbitrage can be simplified.

\begin{prop}\label{prop2.05} Let $Y_t(r)$ be an $n$
parameter yield curve model defined on the open set $O \subset
\R^n$.  Assume that the model satisfies condition SPA with
$\varphi:\R_+ \times O \rightarrow O$ the static price flow
for the model.  Let $x$ be a self-financing investment at $r \in
O$.

The bundle $x$ is a null investment at $r$ if and only if for all $r^\prime$
in some neighborhood of $r$ in $O$
\begin{equation}\label{2.26}
\sum_t \ P_t(r^\prime) x_t \ = 0.
\end{equation}

The bundle $x$ is a local arbitrage at $r$ if and only if for all $r^\prime$
in some neighborhood of $r$ in $O$
\begin{equation}\label{2.27}
\sum_t \ P_t(r^\prime)x_t \ \geq \ 0.
\end{equation}
\end{prop}

{\bf Proof.}  We obtain (\ref{2.26}) and (\ref{2.27}) from (\ref{2.24}) and (\ref{2.25})
respectively by setting $h = 0$.  Thus, (\ref{2.26}) and (\ref{2.27}) are
always necessary.  To prove sufficiency in the SPA case we apply
Lemma \ref{lem2.03} to get $\tilde{\varphi}: \tilde{O} \rightarrow O$, the
extension of the static price flow $\varphi$.  Now let $U_0$ be a
neighborhood of $r$ in $O$ such that $r^\prime \in U_0$ implies
(\ref{2.26}) or (\ref{2.27})  holds.  We can choose
$\epsilon > 0$ and a neighborhood $U$ of $r$ such that for
$(h,r^\prime) \in (-\epsilon,\epsilon) \times U$,
$\tilde{\varphi}_{-h}(r^\prime) \in U_0$.  Shrinking $\epsilon >
0$ we can assume that $x_t = 0$ for $t \leq \epsilon$.  Now for
$(h,r^\prime) \in [0,\epsilon) \times U$ apply (\ref{2.26}) or (\ref{2.27})
with $r^\prime$ replaced by $\tilde{\varphi}_{-h}(r^\prime)$.
Divide through by the common factor $P_h(\tilde{\varphi}_{-
h}(r^\prime))$ and apply (\ref{2.21}) to get (\ref{2.24}) or (\ref{2.25}) .

$\Box$  \vspace{.5cm}

Let $x$ be a bundle of futures.  For an $n$ parameter yield curve
model on $O$ we define the \emph{value function} of $x$, $V: O
\rightarrow \R$ by
\begin{equation}\label{2.28}
V(r) \ = \sum_t \ P_t(r) x_t.
\end{equation}
If $O$ is an open subset of $\R^n$ and the yield curve model
allows static prices then Proposition \ref{prop2.05} says that $x$ is a local
arbitrage at $r$ exactly when $V(r) = 0$ and the value function
$V$ has a local minimum at $r$.  Recall that if $V$ has a local
minimum at $r$ then the partial derivatives $\partial V/ \partial
r_i$ ($i = 1, \ldots, n$) vanish at $r$, i.e. $r$ is a {\em
critical point} for $V$, and the {\em Hessian}, the matrix of
second partials $(\partial^2V/\partial r_i \partial r_j)$, is
positive semidefinite at $r$.  Conversely, if $r$ is a critical
point for $V$ and the Hessian matrix at $r$ is positive definite
then $V$ has a strict local minimum at $r$.
\vspace{.5cm}

\subsection{Linearity (LIN)}

We call a yield curve model \emph{linear}
when the set of admissible yield functions is convex in the
vector space of all smooth real functions on $\R_+$.
Recall that a subset $O$ of a real vector space is {\em convex}
when $r_0,r_1$ in $O$ implies the segment between them is
contained in $O$, i.e. $(1-t)r_0 + tr_1 \in O$ for all $0 \leq t
\leq 1$.  For an $n$ parameter model we want the parametrization
itself to be a linear function.
\vspace{.5cm}

\begin{df}\label{df2.06}  Let $Y_t(r)$ be an $n$
parameter yield curve model defined on $O$.  We say that the
model is {\em linear}, or satisfies condition LIN, if $O$ is a
convex open subset of $\R^n$ and for each $t$ in
$\R_+$ the function $r \mapsto Y_t(r)$ is the restriction
to $O$ of a real-valued linear function from $\R^n$ to $\R$.\end{df}
\vspace{.5cm}

Notice that if $F:O \rightarrow \R$ is the restriction to
$O$, open in $\R^n$, of a linear function, then $F^i =
\partial F/\partial r_i$ is constant on $O$ and $F(r_1, \ldots,r_n) \ = \sum^n_{i=1} \ r_i F^i$.
So if $Y_t(r)$ is a linear $n$
parameter yield curve model on $O$ then $Y^i_t = \partial
Y_t(r)/\partial r_i$ is a smooth real valued function of $t$ which is
independent of $r$ in $O$.  Furthermore, for $r = (r_1, \ldots, r_n)$ in
$O$:
\begin{equation}\label{2.29}
Y_t(r_1, \ldots, r_n) \ = \sum^n_{i=1} \ r_i Y^i_t.
\end{equation}

\begin{lem}\label{lem2.07} For a set $T$ let $\R^T$
denote the vector space of real valued functions defined on $T$.
Let $\{Y^1,\ldots, Y^n\}$ be a list of functions in $\R^T$
for some positive integer $n$ and let $O$ be a nonempty open
subset of $\R^n$.  For $r \in O$ define $Y(r) \ =
\sum^n_{i=1} \ r_i Y^i$ in $\R^T$.  The following conditions are
equivalent
\begin{itemize}
\item[i.] The list $\{Y^1, \ldots, Y^n\}$ is linearly independent in
$\R^T$.

\item[ii.] The function $r \mapsto Y(r)$ from $O$ to $\R^T$ is
injective.

\item[iii.] The set of vectors $\{(Y^1_t, \ldots, Y^n_t): t \in T\}$
spans $\R^n$.
\end{itemize}
\end{lem}

\proof  (i) $\Leftrightarrow$ (ii):  If $Y(r) = Y(r')$
then $\sum_i \ c_i Y^i \ = 0$ where $c_i = r_{i} - r'_{i}$ ($i = 1,
\dots, n)$.  Conversely, if $\sum_i \ c_i  Y^i \ = 0$ and $r \in O$
then because $O$ is open there is an $\epsilon > 0$ small enough
that $r_{\epsilon} \ = r + \epsilon c \in O$.  Clearly, $Y(r) \ =
Y(r_{\epsilon})$.

(i) $\Leftrightarrow $ (iii):  The set of vectors $\{(Y^1_t,
\ldots, Y^n_t)\}$ does not span $\R^n$ if and only if there
is a nonzero vector $c$ in $\R^n$ perpendicular to all of
them, i.e. satisfying $\sum_i \ c_i Y^i_t \ = 0$ for all $t \in T$.

$\Box$ \vspace{.5cm}

{\bf Remark.}  With the same proof we can replace the real field
$\R$ with the complex field $\C$ in the above lemma.
\vspace{.5cm}

\begin{df}\label{def2.07a}For a linearly independent list $\{Y^i: i = 1, \ldots, n\}$ of
smooth functions in $\R^{\R_+}$  the \emph{associated linear model} is the linear map
$Y: \R^n \rightarrow \R^{\R_+}$ given by (\ref{2.29}), i.e.
$Y(r) \ = \sum_i \ r_i Y^i$.

The associated {\em positive domain} and
{\em strict positive domain} are the subsets $O$ and $O_+$
defined by
\[O(Y^1, \ldots, Y^n) \ = \{r \in \R^n: Y^i_t(r) > 0 \; {\rm for
\; all\;}t \in \R_+, i = 1,\ldots,n \}\]
\begin{equation}\label{2.30}
O_+(Y^1,\ldots, Y^n) \ = \{r \in \R^n: \inf_t Y^i_t(r) > 0, i = 1,\ldots,n \}.
\end{equation}
\end{df}
\vspace{.5cm}

So $r \in O_+(Y^1, \ldots, Y^n)$ exactly when for some $\delta > 0$,
 $Y^i_t(r) \geq \delta$ for all $ t \in \R_+, i = 1,\dots,n$.  It easily
follows that $O(Y^1,\ldots, Y^n)$ and $O_+(Y^1, \ldots,
Y^n)$ are convex subsets of $\R^n$.
\vspace{.5cm}

\begin{lem}\label{lem2.08}If $\{Y^i: i = 1, \ldots, n\}$ is
a linearly independent list of bounded, smooth functions in $\R^{\R_+}$,
then $O_+(Y^1, \ldots, Y^n)$ is open as well as
convex.\end{lem}

\proof Assume $|Y^i_t| \leq M$ for $i = 1, \ldots, n$ and
$t \in \R_+$.  If $Y^i_t(r) \geq \delta$ for all $t$ with
$\delta > 0$ and $|r'_i - r_i| \ < \ \delta/2n M$ for $i = 1,
\ldots, n$, then $Y^i_t(r') \geq \delta/2$ for all $t$.

$\Box$ \vspace{.5cm}

\begin{prop}\label{prop2.09}If $Y_t(r)$ is a linear $n$
parameter yield curve model defined on a convex open set $O$,
then the list $\{Y^1, \ldots, Y^n\}$ of smooth functions in $\R^{\R_+}$ defined by
\begin{equation}\label{2.31}
Y^i_t \ = \partial Y_t(r)/\partial r_i \;\;\; (i = 1, \ldots, n)
\end{equation}
is linearly independent.

Conversely, if $\{Y^1, \ldots, Y^n\}$ is a linearly independent
list of smooth functions in $\R^{\R_+}$ then the associated linear model
 defines a linear yield curve model on any nonempty,
convex, open set $O$ contained in the positive domain $O(Y^1, \ldots, Y^n)$.\end{prop}

\proof By Lemma \ref{lem3.07} the map $Y: O \rightarrow \R^{\R_+}$ is injective iff $\{Y^1, \ldots, Y^n\}$ is
linearly independent.  The results then follow from Definition 1.

$\Box$ \vspace{.5cm}

{\bf Remark.}  It follows that for any linear yield curve model
the function $Y: O \rightarrow \R^{\R_+}$ extends
uniquely as a linear map from $\R^n$ to ${\bf
R}^{\R_+}$.  Thus, we use the formula (\ref{2.29}) to define
$Y(r)$ for all $r$ in $\R^n$.  In particular, a linearly independent list
$\{Y^1,\ldots,Y^n\}$ defines a linear yield curve model for some
$O$ exactly when the interior of $O(Y^1,\ldots Y^n)$ in $\R^n$ is nonempty.
\vspace{.5cm}

To combine conditions SPA and LIN we recall that for an $n \times
n$ matrix $A = (A^i_j)$ the smooth matrix function defined for $h
\in \R$:
\begin{equation}\label{2.32}
E_h \ = \exp[hA] \ = \sum^\infty_{N=0} \ (hA)^N/N!
\end{equation}
is the {\em fundamental solution} for the linear system
associated with $A$.  That is, differentiating with respect to $h$:

\[E_0 \ = I,\;\; {\rm the \; identity \; matrix}\]
\[\stackrel{\bullet}{E}_h  \ = E_hA \ = AE_h,\; {\rm i.e.}\]
\begin{equation}\label{2.33}
\stackrel{\bullet}{E^i_{jh}}\ = \sum^n_{k=1} \ E^i_{kh} A^k_j \ =
\sum^n_{k=1} \ A^i_k E^k_{jh}.
\end{equation}
For $r \in \R^n$,  define $\psi_h(r) \ = rE_h$ so that
$\psi_h(r)_j \ = \sum^n_{i=1} \ r_i E^i_{jh}$.  Then $\psi_h(r)$ is
the solution of the linear initial value problem on $\R^n$.
\[\psi_0(r) \ = r\]
\begin{equation}\label{2.34}
\stackrel{\bullet}{\psi}_h(r) \ = \psi_h(r)A
\end{equation}
That is, regarded as a function from $\R \times \R^n$
to $\R^n$,  $\psi$ is the linear solution flow for the linear
vectorfield defined by $\xi(r) \ = rA$.

\vspace{.5cm}

\begin{prop}\label{prop2.10} Suppose that $\varphi:
\R_+ \times O \rightarrow O$ is the static price flow of a
linear $n$ parameter yield curve which is defined on the open
convex set $O$ and which satisfies condition SPA.  For each $h
\geq 0$, $\varphi_h$ is the restriction to $O$ of some linear
automorphism of $\R^n$.  In fact, there is an $n \times n$
matrix $A$ such that for $(h,r) \in \R_+ \times O$,
$\varphi_h(r) = r E_h$ where $E$ is the fundamental solution for
the matrix differential equation associated with $A$.\end{prop}

\proof  By hypothesis the map $Y: O \rightarrow \R^{\R_+}$
associating to $r$ the function $Y(r)$ is the
restriction of a linear map $\tilde{Y}: \R^n \rightarrow
\R^{\R_+}$, and $\tilde{Y}$ is injective because its
restriction to $O$ is (apply Lemma \ref{lem2.07}).  For $h \geq 0$ the
translation map $\Phi_h: \R^{\R_+} \rightarrow \R^{\R_+}$
defined by (\ref{2.14}) is clearly linear.  Because
$\Phi_h \circ Y \ = Y \circ \varphi_h$ on $O$ by (\ref{2.15}) we can
define the linear map $\tilde{Y}^{-1} \circ \Phi_h \circ
\tilde{Y}: \R^n \rightarrow \R^n$ which restricts to
$\varphi_h$ on $O$.

Because $\varphi_h$ is linear in the $r$ variables for each $h$,
the same is true of the derivative, the vectorfield $\xi: O
\rightarrow \R^n$ defined by (\ref{2.19}).  That is, there exists
an $n \times n$ matrix $A$ such that
\begin{equation}\label{2.35}
\xi(r) \ = rA.
\end{equation}
Comparing (\ref{2.34}) with (\ref{2.20}), we see that uniqueness for
solutions of smooth differential equations implies that
$\varphi_h$ is the restriction to $O$ of $\psi_h.$

$ \Box$\vspace{.5cm}

\begin{cor}\label{cor2.11} Let $Y_t(r)$ be a linear $n$
parameter yield curve model defined on $O$ which satisfies
condition SPA.  Let $A$ be the $n \times n$ matrix whose
fundamental solution $E$ defines the static price flow as
described in Proposition \ref{prop2.10}.  There is a vector $(Y^1_0, \ldots,Y^n_0)$
 in $\R^n$ such that the short rate function $Y_0(r)$
is given by $\sum^n_{i=1} \ r_i Y^i_0$ and for all $(t,r) \in
\R_+ \times O$, and the basic yield curves satisfy:

\[Y^i_t \ = \sum^n_{k=1} \ E^i_{kt} Y^k_0\]
\begin{equation}\label{2.36}
\stackrel{\bullet}{Y}_t^i \ = \sum^n_{k=1} \ A^i_k Y^k_t.
\end{equation}\end{cor}

\proof For the short rate apply (\ref{2.29}) with $t = 0$.  For
the first equation of (\ref{2.36}) apply (\ref{2.18}) with $\varphi(r) \ = rE_t$.
Then use the second equation in (\ref{2.33}) to compute
$\stackrel{\bullet}{Y}^i_t$.

$\Box$\vspace{.5cm}
%

\subsection{Long Rates Exist (LRE)}

For any yield curve $Y_t$ the \emph{short rate} is the positive number $Y_0$.
We define the \emph{long rate} to be $Lim_{t \rightarrow \infty} \ Y_t$ if this limit
exists as a finite real number.  If the long rate exists then $Y$
is bounded on $\R_+$.  A bounded function may however fail
to have a limit $t$ as $t$ approaches infinity.  The function may
instead oscillate.

\vspace{.5cm}

\begin{df}\label{df2.12} Let $Y_t(r)$ be an $n$
parameter yield curve model defined on $O$.  We say that the
model has {\em long rates}, or satisfies condition LRE, if for
each $r \in O$
\begin{equation}\label{2.37}
Y_\infty(r) \ \equiv \ Lim_{t \rightarrow \infty}\ \ Y_t(r)
\end{equation}
exists and is finite.\end{df}

\vspace{.5cm}

{\bf Remark.}  If the model satisfies LRE then for each $r \in
O$, $Y(r)$ is a bounded smooth function in $\R^{\R_+}.$
\vspace{.5cm}

As usual, for the linear case there are special results.
\vspace{.5cm}

\begin{prop}\label{prop2.13} Let $Y_t(r) \ = \sum_i \ r_i Y^i_t$ be a
linear $n$ parameter yield curve model.  The model
has long rates if and only if for $i = 1, \ldots, n$, the limit
$Y^i_{\infty} \ = Lim_{t \rightarrow \infty} \ Y^i_t \ $ exists.  In
that case the function $Y_\infty(r)$ is well defined by (\ref{2.37})
for all $r \in \R^n$ and
\begin{equation}\label{2.38}
Y_{\infty}(r) \ = \sum^n_{i=1} \ r_i Y^i_{\infty}.
\end{equation}\end{prop}

\proof Those $Y$ in $\R^{\R_+}$ for which the
limit $Y_\infty$ exists constitute a linear subspace.  This
subspace contains the image of the linear map from $\R^n$ to
$\R^{\R_+}$ given by the model when it contains the
$Y(r)$'s as $r$ varies over some basis for $\R^n$.  So
containing the $Y^i$'s is sufficient as is containing the
$Y(r)$'s for $r$ varying over some open set in $\R^n$.
Equation (\ref{2.38}) is then clear from the linearity of taking
limits.

$\Box$ \vspace{.5cm}


\subsection{Main Results}

In comparing two yield curve models we say that the first is
included in the second if every yield curve admissible with
respect to the first model is admissible for the second.  In
comparing two linear models we have special results.
\vspace{.5cm}

\begin{df}\label{df2.14} Assume that $\{ Y^1, \ldots, Y^n \}$ and $\{ \tilde{Y}^1, \ldots, \tilde{Y}^m \}$
are linearly independent lists of smooth functions in $\R^{\R_+}$ with associated linear models $Y_t(r)  = \sum^n_{i=1} \ r_iY^i_t$ and
$\tilde{Y}_t(s)  = \sum^m_{j=1} \ s_j \tilde{Y}^j_t$.  We say that the
$\tilde{Y}_t(s)$ model {\em includes} the $Y_t(r)$ model when the following equivalent conditions hold.
\newpage
\begin{enumerate}
\item[(1)]  There exists a nonempty open subset $O_1 \subset \R^n$
such that for each $r \in O_1$ there exists $s \in \R^m$
such that $Y_t(r) \ = \tilde{Y}_t(s)$ for all $t \in \R_+$.

\item[(2)]  There is an injective linear mapping $S: \R^n \rightarrow \R^m$
such that for all $r \in \R^n$, $Y(r)
= \tilde{Y}(S(r))$.
\end{enumerate}

In that case, $S$ maps $O(Y^1,\ldots, Y^n)$ into $O(\tilde{Y}^1,\ldots, \tilde{Y}^m)$
and $O_+(Y^1,\ldots, Y^n)$ into $O_+(\tilde{Y}^1,\ldots, \tilde{Y}^m)$.
\end{df}
\vspace{.5cm}

{\bf Proof of equivalence.}  (1) $\Leftrightarrow$ (2):  By
assumption the linear maps $Y: \R^n \rightarrow \R^{\R_+}$ and
$\tilde{Y}: \R^m \rightarrow \R^{\R_+}$ are injective.  Assumption (1) says that on the
nonempty open set $O_1$ the map $\tilde{Y}^{-1} \circ Y$ is well-defined.  
As $O_1$ spans $\R^n$ the image of $Y$ must lie in
the image of $\tilde{Y}$ by linearity.  Hence, $\tilde{Y}^{-1}
\circ Y: O_1 \rightarrow \R^m$ extends to define the linear
map $S$.

The converse is obvious.

Since $Y(r) = \tilde{Y}(S(r))$ it is clear that $r$ is in $O$ or $O_+$ if and only if
$S(r)$ is.

$\Box$\vspace{.5cm}
%

We are now ready to state our main results.
\vspace{.5cm}

\begin{theo}\label{theo2.15} Let $Y_t(r)$ be an $n$
parameter yield curve model which is linear and allows static
prices, i.e. $Y_t(r)$ satisfies SPA and $LIN$.  If the model
allows no local arbitrage, i.e. satisfies condition NLA, then it
includes one of the following four  linear models:
\begin{enumerate}

\item[(1)]  With $\rho \neq 0$ fixed, let
\begin{equation}\label{2.39}
Y_t(r_1,r_2) \ = r_1 e^{\rho t} + r_2 e^{2\rho t}.
\end{equation}

\item[(2)]  With $\rho \neq 0$ and $\omega > 0$ fixed, let
\begin{equation}\label{2.40}
Y_t(r_1,r_2,r_3) \ = r_1 e^{\rho t} \cos \omega t + r_2 e^{\rho t}
\sin \omega t + r_3 e^{2 \rho t}.
\end{equation}

\item[(3)]  Let
\begin{equation}\label{2.41}
Y_t(r_1,r_2) \ = r_1 + r_2t.
\end{equation}


\item[(4)]  With  $\omega > 0$ fixed, let
\begin{equation}\label{2.42}
Y_t(r_1,r_2,r_3) \ = r_1 \cos \omega t + r_2 \sin \omega t.
\end{equation}
\end{enumerate}

Conversely, if the yield curve model includes either model (1), (2) or (3)
then it satisfies NLA. \end{theo}
\vspace{.5cm}

The proof of all this occupies the next section.  We will also
describe the oddities associated with model (4).  Because model
(3) contains unbounded functions and model (4) has undamped
oscillation we regard them as anomalous.  Also in (1) and (2), $\rho > 0$
leads to unbounded functions.

There are some positivity problems in the above linear models.  In  (1) if $\rho$ is negative, the set $O_+(e^{\rho t},e^{2\rho t})$ is empty,
although $\{ (r_1,r_2) : r_1, r_2 > 0 \}$ is contained in
$O(e^{\rho t},e^{2\rho t})$. In (2) we can write $(r_1,r_2)$ in polar coordinates as $(a \cos \phi,a \sin \phi )$ and then model (2) can be written as
$ e^{\rho t}[ a \cos (\omega t - \phi) +  r_3 e^{\rho t}]$. It follows that if $\rho$ is negative, then
$O(e^{\rho t} \cos \omega t,e^{\rho t}\sin \omega t,e^{2 \rho t})$ is empty.  Similarly, in (4) $O( \cos \omega t,\sin \omega t)$
is empty. These problems are eliminated if we also demand that
the model includes the flat rate model.  That is,
\begin{align}\begin{split}&O_+(1, e^{\rho t},e^{2\rho t}),\\
&O_+(1, e^{\rho t} \cos \omega t,e^{\rho t}\sin \omega t,e^{2 \rho t}),\\
&O_+(1,t)\quad \text{and} \quad
O_+(1, \cos \omega t,\sin \omega t)
\end{split}\end{align}
are nonempty for all $\r$ and $\omega$.

The two simple yield
curve models in the title come from the following:
\vspace{.5cm}

\begin{cor}\label{cor2.16} Let $Y_t(r)$ be an $n$
parameter yield curve model which is linear, allows static prices
and has long rates, i.e. $Y_t(r)$ satisfies SPA, LIN and LRE.
The model contains the flat yield curve model and allows no local
arbitrage, i.e. satisfies NLA, if and only if it includes one of
the following two yield curve models.
\begin{enumerate}

\item[(1)]  {\bf Exponential Model} (3 parameters):  With $\rho > 0$ fixed, let
\begin{equation}\label{2.43}
Y_t(r_1,r_2,r_3) \ = r_1 + r_2 e^{-\rho t} + r_3 e^{-2\rho t},
\end{equation}
defined on the open convex set $O_+(1,e^{-\rho t},e^{-2 \rho
t})$.  The short rate is $r_1 + r_2 + r_3 $ and the long rate is
$r_1$.

\item[(2)]  {\bf Exponential-Oscillation Model} (4 parameters):  With
$\rho, \omega > 0$ fixed, let
\[Y_t(r_1,r_2,r_3,r_4) \ = r_1 + r_2 e^{-\rho t} \cos \omega t +\]
\begin{equation}\label{2.44}
r_3 e^{-\rho t} \sin \omega t + r_4 e^{-2 \rho t},
\end{equation}
defined on the open convex set $O_+(1,e^{-\rho t} \cos \omega t,e^{-\rho t}\sin \omega t,e^{-2 \rho t})$
The short rate is $r_1 + r_2 + r_4$ and the long rate is $r_1$.
\end{enumerate}

By choosing time units property, i.e. by a constant rescaling, we
can assume $\rho = 1$ in the models.\end{cor}

\proof  In each case the $r_1$ coordinate comes from the
assumption that the flat yield curve model is included.  Theorem \ref{theo2.15}
 says that any model which includes one of these satisfies NLA
because (1) and (2) here contain (1) and (2) of the theorem,
respectively.  Conversely, condition LRE excludes unbounded
growth, i.e. model (3) of the theorem and $\rho > 0$ in models
(1)and (2).  It also excludes the undamped oscillation of model (4).
So the theorem implies either model (1) or (2) here is included.

The rescaling in the last paragraph is the replacement of the
time variable $t$ by $\rho t$.  $\hfill \Box$
\vspace{1cm}

\section{Arbitrage in Linear Models}
\vspace{.5cm}

Linear models occur as subspaces of the real vector space $\R^T$ for some set $T$.  As $\R$ is a subspace of $\C$, we can regard
$\R^T$ as a real subspace of the complex vector space $\C^T$.  For $g \in \C^T$ we can write $g = f_1 + \ii f_2$, with $f_1, f_2 \in \R^T$ and we
write $f_1 = Re(g), f_2 = Im(g)$ with the conjugate $\bar g = f_1 - \ii f_2$.

For $V$ a subspace of $\R^T$ we define $V_{\C}$ to be the complex subspace
of $\C^T$ given by:
\begin{equation}\label{eq3.001}
V_{\C} \ = \ \{ f_1 + \ii f_2: f_1, f_2 \in V \}.
\end{equation}
So, for example, $\C^T = (\R^T)_{\C}$.

\begin{lem}\label{lem3.001} Let $W$ be a complex subspace of $\C^T$.  There exists a real subspace $V$ of $\R^T$ such that $W = V_{\C}$ if and
only if $W$ satisfies the following equivalent conditions.
\begin{itemize}
\item[(i)]  If $g \in W$, then $\bar g \in W$.
\item[(ii)]  If $g \in W$, then $Re(g) \in W$.
 \item[(iii)]  If $g \in W$, then $Re(g), Im(g) \in W$.
 \end{itemize}

 In that case,
 \begin{equation}\label{eq3.002}
 V \ = \ W \cap \R^T \ = \ \{ Re(g) : g \in W \}  \ = \ \{ Re(g), Im(g) : g \in W \},
 \end{equation}
 and for any index set $I$, if $\{ g_i : i \in I \}$ spans $W$, then $\{ Re(g_i), Im(g_i) : i \in I \}$ spans $V$.
 \end{lem}
 \vspace{.5cm}

{\bf Proof:} (i) $\Rightarrow$ (ii): $Re(g) = \frac{1}{2}(g + \bar g)$.

 (ii) $\Rightarrow$ (iii): $Im(g) = - Re(\ii g)$.

 (iii) $\Rightarrow$ (i): $\bar g = Re(g) - \ii Im(g).$

 Clearly, $V_{\C}$ satisfies (iii). On the other hand, given (i)-(iii) it is clear that $W = V_{\C}$ for $V$ given by (\ref{eq3.002}).

 Finally, $g \ = \ \sum_i \ a_i g_i$ implies $Re(g) = \sum_i Re(a_i) Re(g_i) - Im(a_i) Im(g_i)$ (with $a_i = 0$ except for finitely many $i \in I$).

 $\Box$ \vspace{.5cm}

Throughout this section we consider linear models of the form \\
$Y_t(r_1, \ldots, r_n) \ = \sum^n_{i=1} r_iY^i_t$, which we regard
as a linear mapping $Y: \R^n \rightarrow \R^{\R_+}$ with  $\{ Y^1_t(r), \dots, Y^n_t(r) \}$ a linearly independent list of
smooth functions in $\R^{\R_+}.$
For a yield curve model, we restrict to some nonempty convex open set $O \subset \R^n$
such that $r \in O$ implies
$Y^i_t(r) > 0$ for all $t \geq 0$ and $i = 1, \dots, n$.  We can regard the log-price as
the linear function $L: \R^n \rightarrow \R^{\R_+}$ given by
\[L_t(r_1, \ldots, r_n) \ = \sum^n_{i=1} \ r_i L^i_t\]
\begin{equation}\label{3.01}
L^i_t \ = \int^t_0 \ Y^i_s \ ds.
\end{equation}
So, in particular, $L_0(r) = 0$ for all $r \in \R^n$.  The
price-function extends as well to a nonlinear function $P: \R^n \rightarrow \R^{\R_+}$ by
\begin{equation}\label{3.02}
P_t(r_1, \ldots, r_n) \ = \exp[-\sum^n_{i=1} \ r_i L^i_t].
\end{equation}
For $r \in O$ the price and log-price curves satisfy (\ref{2.04}) and
(\ref{2.07}).

Because of the linear structure of the models we can describe
investment procedures which might lead to an arbitrage
possibility at every $r \in O$.

Let $z_t$ be a finite function on $(0,\infty)$.  That is, it is a
real-valued function on $\R_+$ whose support is a finite
subset excluding $0$.  Since the function $z_t$, once chosen,
will be fixed during some extended arguments it will be useful to
introduce the notation:
\begin{equation}\label{3.03}
\langle f \rangle \ = \sum_t \ z_t f(t)
\end{equation}
for any function $f$ in $\R^{\R_+}$ or in $\C^{\R_+}$.

Think of $z_t$ as the positive or negative amount of current
dollars which we allocate to futures $t$ time units from now.  If
the current price curve is $P_t(r^*)$ then the quantity $x_t$ of
time $t$ futures with current worth $z_t$ is given by
\begin{equation}\label{3.04}
x_t \ = z_t/P_t(r^*).
\end{equation}
So the present value of the bundle of futures $x_t$ is given by
\begin{equation}\label{3.05}
\sum \ x_t P_t(r^*) \ = \langle 1 \rangle.
\end{equation}
For any $r \in \R^n$ the value function, defined by (\ref{2.28}),
is at $r$ given by
\begin{equation}\label{3.06}
V(r) \ = \langle P(r-r^*)\rangle
\end{equation}
because for each $t$ $P_t(r) = P_t(r^*) P_t(r-r^*)$ by (\ref{3.02}).
Taking the partial derivative with respect to $r_i$ we see that
\begin{equation}\label{3.07}
\frac{\partial V}{\partial r_i}|_{r = r^*} \ = - \langle L^i \rangle
\end{equation}
and the Hessian matrix is given by
\begin{equation}\label{3.08}
\frac{\partial V}{\partial r_i \partial r_j}|_{r = r^*} \ = \langle L^iL^j \rangle.
\end{equation}

Thus, when $r = r^*$ the investment $x_t$ is self-financing provided $\langle 1 \rangle = 0$.  
It is then a critical point for the value function
when $\langle L^i \rangle = 0$ for $i = 1, \ldots, n$.  It is a
strict local minimum when the $n \times n$ Hessian matrix
$\langle L^iL^j\rangle$ is positive definite.  By Proposition
\ref{prop2.05}, if the linear yield curve  model satisfies SPA, then for $r^* \in O$ such a local minimum
is a local arbitrage for the model.

We can reverse this analysis as well.  If $x_t$ is a bundle of
futures with current prices given by $P_t(r^*)$ we can use (\ref{3.04})
to define the finite function $z_t$.

The following \emph{No Arbitrage Lemma} does not require the SPA
condition.  We postpone the proof.
\vspace{.5cm}

\begin{lem}\label{lem3.01} For the $n$ parameter linear yield
curve model whose log-price function $L$ is given by
(\ref{3.01}) let $L(\R^n)$ be the linear subspace of $\R^{\R_+}$
which is the image of $L$.  Assume that for some
positive integer $p$ there exist functions $F^1, \ldots, F^p$ in
$L(\R^n)$ which satisfy the following conditions:
\begin{enumerate}

\item[(1)]  The function $t \mapsto (F^1_t, \ldots, F^p_t)$ is an
injective function from $\R_+$ to $R^p$, i.e. if $t_1
\neq t_2$ then for some $j = 1,\ldots, p$, $F^j_{t_1} \neq F^j_{t_2}$.

\item[(2)]  The function $\sum^p_{j=1} \ (F^j)^2$, the sum of the squares,
is in the subspace $L(\R^n)$.

\end{enumerate}

The model satisfies condition NLA.  In fact, if $z_t$ is a finite
function on $(0,\infty)$ not identically zero but satisfying
$\langle 1 \rangle = 0$ and $\langle L^i\rangle = 0$ for $i =
1,\ldots, n$, then for any $r^*$ in $\R^n$, the value
function defined by using $x_t = z_t/P_t(r^*)$ has a critical
point at $r = r^*$ which is not a local minimum.\end{lem}
\vspace{.5cm}

Now we assume that the linear model satisfies condition SPA.  By
Corollary \ref{cor2.11} the yield functions are solutions of linear
ordinary differential equations with constant coefficients.  For
such functions it is useful to introduce some special notation.

In the additive group $\Z \times \C$ the set of \emph{exponents} $\E$
consists of the elements whose integer
coordinate is nonnegative.  So an exponent $q$ is a pair
$(m,\lambda)$ with $m$ a nonnegative integer and $\lambda$ a
complex number.  For $q = (m,\lambda)$ with $\lambda = a + ib$ we
define the {\em conjugate} $\overline{q} =
(m,\overline{\lambda})$, the {\em real part} $\rho(q) = (m,a)$
and the {\em imaginary part} $\omega(q) = b$.  Define
\begin{equation}\label{3.09}
q^\prime \ = (\max(0,m-1),\lambda) \quad \mbox{for} \quad q \ = (m,\lambda).
\end{equation}
so that $q^\prime = q$ when $q = (0,\lambda)$.

For the exponent
$q = (m,\lambda)$ define the function $f^q: \R \rightarrow \C$ by
\[f^q(t) = t^m \exp (\lambda t)\]
\begin{equation}\label{3.10}
\ =  f^{\rho(q)}(t)[\cos(\omega(q)t) + i \sin(\omega(q)t)].
\end{equation}
Clearly $f^{\overline{q}} \ = \overline{f^q}$ and $f^{q_1 + q_2} \ = f^{q_1} \cdot f^{q_2}$
for $q,q_1,q_2 \in \E$.  In
addition, the derivative satisfies
\begin{equation}\label{3.11}
\stackrel{\bullet}{f^q} \ = \lambda f^q + mf^{q^\prime} \;\;\; (q =
(m,\lambda)).
\end{equation}

A \emph{closed set} of exponents $Q$ is a finite subset of $\E$ such
that $q \in Q$ implies $\overline{q}$, $q^\prime \in Q$.
For a closed set $Q$ and a complex number $\lambda$ we define the
\emph{multiplicity} of $\lambda$ in $Q$ to be
\begin{equation}\label{3.12}
\min\{m \geq 0: (m,\lambda) \not\in Q\}.
\end{equation}
\vspace{.5cm}

\begin{lem}\label{lem3.02} (a)  Let $Q$ be any finite subset
of $\E$ and $T$ be any nonempty open subset of $\R$.
The set of restrictions to $T$ of the functions $\{ f^q: q \in Q \}$ is
linearly independent in the complex vector space $\C^T$.

(b)  Let $V$ be a finite dimensional subspace of $\R^{\R_+}$
consisting of differentiable functions.  If $V$
is spanned by $\{{\rm Re}(f^q),{\rm Im}(f^q): q \in Q\}$ for some
closed set of exponents $Q$, then $V$ is closed under
differentiation, i.e. $f \in V$ implies $\stackrel{\bullet}{f} \in V$.
Conversely, if $V$ is closed under differentiation, then
$Q \ = \{ q \in \E: { \rm Re}(f^q),{\rm Im}(f^q) \in V\}$ is a
closed set of exponents such that $\{{\rm Re}(f^q),{\rm Im}(f^q):
q \in Q\}$ spans $V$.  In particular, if $\{ Y^1,\ldots, Y^n \}$ is
a list of differentiable functions such that
\begin{equation}\label{3.13}
\stackrel{\bullet}{Y^i_t} \ = \sum^n_{k=1} \ A^i_k Y^k_t \;\;\;(i =
1,\ldots, n),
\end{equation}
then the vector space $V$ spanned by $\{Y^1, \ldots, Y^n\}$ is
spanned by the real and imaginary parts of $f^q$ for $q$ varying
in such a closed set of exponents.\end{lem}

\proof (a)  Denote by $D$ the derivative operator on the
subspace of differentiable functions in $\C^{\R_+}$,
i.e. $Df = \stackrel{\bullet}{f}$.  Suppose that $\sum \ c_q   f^q$
is identically zero on the open set $T$.  Let $q_0 =
(m_0,\lambda_0)$ be an exponent
such that $c_{q_0} \neq 0$.  For each $\lambda \neq \lambda_0$
apply the operator $(D - \lambda)$ $m$ times where $m$ is the
multiplicity of $\lambda$ in $Q$, and $D - \lambda_0$ $m_0 - 1$
times.  The only term remaining in the relation is $c_{q_0}$ times a
nonzero number times $f^{(0,\lambda_0)}$. This term is not identically
zero on $T$ unless $c_{q_0} = 0$ after all.\vspace{.25cm}

(b)  The first statement is obvious by (\ref{3.11}).  For the converse,
assume $V$ is closed under differentiation so that  the complex
subspace $V_{\C} = \{g \in \C^{\R_+}: {\rm
Re}(g),{\rm Im}(g) \in V\}$ and  $Q \ = \{ q \in \E: f^q \in V_{\C} \}$.
 Since $V$ and hence $V_{\C}$ are finite
dimensional $Q$ is a finite subset of $\E$ by (a).  If $q
\in Q$, then $\overline{q} \in Q$ because $V_{\C}$ is closed
under conjugation.  By (\ref{3.11}) $q^\prime \in Q$ when $m > 0$. Since $q^\prime =
q$ when $m = 0$, it follows that $Q$ is closed.

The differentiation operator $D$ is a linear map of $V_{\C}$
to itself and so we can choose a basis $\{Y^1 , \ldots, Y^n \}$ so
that $D$ is represented by a matrix in Jordan canonical form or,
to be more precise, so that each $k \times k$ Jordan bloc is of
the form
\begin{equation}\label{3.14}
\left ( \begin{array}{ccccc}
\lambda & k-1 & 0 & \ldots & 0\\
0 & \lambda & k-2 & \ldots & 0\\
0 & \ldots & \ldots & \lambda & 1\\
0 & \ldots & \ldots & \ldots & \lambda
\end{array} \right).
\end{equation}
So we see from (\ref{3.11}) that each $Y^i$ can be chosen a function
$f^q$ with $q = (m,\lambda)$ and $\lambda$ an eigenvalue of the
matrix.  Consequently, the $f^q$'s span $V_{\C}$ and so the
Re$(f^q)$'s and Im$(f^q)$'s span $V$ by Lemma \ref{lem3.001}.

$\Box$\vspace{.5cm}

{\bf Remark.}  Regard $q \in \C^n$ and $t \in \R^n$ as
column vectors i.e. $n \times 1$ matrices, and define
\begin{equation}\label{3.15}
F^q(t) \ = \exp[q^T t] \ = \exp[\sum^n_{i=1} \ q_i t_i].
\end{equation}
For any finite subset $Q$ of $\C^n$ and any nonempty open
subset $O$ of $\R^n$ the set of restrictions to $O$ of the
functions $\{ F^q: q \in Q \}$ is linearly independent in the
complex vector space $\C^O$.  The proof proceeds as in (a)
using partial derivative operators.
\vspace{.5cm}

\begin{cor}\label{cor3.03} For the $n$ parameter linear
yield curve model whose  yield function is $Y: \R^n
\rightarrow \R^{\R_+}$, let $Y(\R^n)$ denote the
linear subspace of $\R^{\R_+}$ which is the image of
$Y$.  If the model satisfies condition SPA, then there is a closed
set of exponents $Q$ such that
\[\{f^q: q \in Q \; {\rm with} \; \omega(q) = 0\} \cup\]
\begin{equation}\label{3.16}
\{{\rm Re}(f^q),{\rm Im}(f^q): q \in Q \; {\rm with} \;
\omega(q) > 0\}
\end{equation}
is a basis for $Y(\R^n)$.\end{cor}

\proof The yield curve function is given by
(\ref{2.29}) with $\{Y^1,\ldots, Y^n\}$ linearly independent by
Proposition \ref{prop2.09}.  So letting $V = Y(\R^n)$ we see that $V$
is an $n$ dimensional subspace of functions to which part (b) of
the lemma applies because Corollary \ref{cor2.11} implies (\ref{3.13}).  It
follows that there is a closed set of exponents $Q$ such that the
functions listed in (\ref{3.16}) span $Y(\R^n)$.  The functions
$\{f^q: q \in Q\}$ span the complex space $V_{\C}$ and so by
part (a) of the lemma they form a basis for $V_{\C}$.  If
$\omega(q) = 0$ then $f^q$ is a real function and if $\omega(q) >
0$ then Re$(f^q) \ = (f^q + f^{\overline{q}})/2$ while Im$(f^q) \ = (f^q - f^{\overline{q}})/2i$.
So the real functions listed in (\ref{3.16})
form a linearly independent set as well.

$\Box$\vspace{.5cm}

The set $Q$ of exponents in Corollary \ref{cor3.03} is given by:
\begin{equation}\label{3.17}
Q \ = \{q \in \E: {\rm Re}(f^q),{\rm Im}(f^q) \in Y(\R^n)\}.
\end{equation}
We call this set the {\em exponent set} for the linear yield
curve model $Y_t(r)$, defined when the model satisfies SPA.

Now for any exponent $q$ define the exponent $\tilde{q}$ by
\begin{equation}\label{3.18}
\tilde{q} \ = \left \{ \begin{array}{cll}
q & {\rm if}\; q = (m,\lambda) & {\rm with} \; \lambda \neq 0\\
(m+1,0) & {\rm if} \; q = (m,0) &
\end{array} \right.
\end{equation}
In particular, $\tilde{0} = (1,0)$.  For any set of exponents $Q$
let $\tilde{Q} \ = $ \\ $ \{\tilde{q}: q \in Q\}$.  Clearly, $0 \not\in
\tilde{Q}$.

For $q$ any nonzero exponent, i.e. $q \neq (0,0)$, define
\begin{equation}\label{3.19}
l^q(t) \ = f^q(t) - f^q(0),
\end{equation}
so that $l^q = f^q$ if $m > 0$ and $l^q = f^q - 1$ if $m = 0$.

\vspace{2ex}

\begin{lem}\label{lem3.04}  Assume that $Y: \R^n \rightarrow \R^{\R_+}$
and $L: \R^n \rightarrow \R^{\R_+}$ are the  yield and log-price
functions for a linear yield curve model satisfying SPA.  Let $Q$ be the
exponent set for the yield curve model.

\[\{l^q: q \in \tilde{Q} \; {\rm with} \; \omega(q) = 0\} \cup\]
\begin{equation}\label{3.20}
\{{\rm Re}(l^q),{\rm Im}(l^q): q \in \tilde{Q} \; {\rm with} \;
\omega(q) > 0\}
\end{equation}
is a basis for $L(\R^n)$.\end{lem}

\proof From (\ref{3.11}) we have, for $q = (m,\lambda)$:
\begin{equation}\label{3.21}
l^q(t) \ = \lambda \int^t_0 f^q(s)ds + m \int^t_0
f^{q^\prime}(s)ds.
\end{equation}
Using induction on $m$, it easily follows that $\{l^q: q \in \tilde{Q}\}$ is a basis for
the complex vector space $V_{\C}$ where $V = L(\R^n)$.
So the list in (\ref{3.20}) is a basis for the real vector space
$L(\R^n)$ just as in Corollary \ref{cor3.03}.

$\Box$ \vspace{.5cm}

We will call $\tilde{Q}$ the {\em adjusted exponent set} for the
yield curve model when $Q$ is the exponent set.

We are now ready to state our second technical result, the {\em
Arbitrage Everywhere Lemma.}

\vspace{2ex}

\begin{lem}\label{lem3.05}  For a linear $n$ parameter yield
curve model defined on $O$ and which satisfies SPA, let
$\tilde{Q}$ be the adjusted exponent set.  Assume that for all $q
\in \tilde{Q}$, $q + \overline{q} \not\in \tilde{Q} \cup \{0\}$.
The model admits strict local arbitrage at every parameter value
$r \in O$.  In fact, there is a finite function $z_t$ whose
support is a subset of $(0,\infty)$ consisting of at most
$\frac{1}{2}(n+1)(n+2)$ times, so that for any $r^*$ in $\R^n$
 the value function defined using $x_t \ = z_t / P_t(r^*)$ has
a strict local minimum at $r = r^*$.\end{lem}

$\Box$ \vspace{.5cm}

Again we postpone the proof, first showing how Lemmas \ref{lem3.01} and \ref{lem3.05}
yield the main result, stated in the previous section.

\vspace{.5cm}

{\bf Proof of Theorem \ref{theo2.15}:}  Since the model is linear and
satisfies SPA we can apply the contrapositive of Lemma \ref{lem3.05}.  When
the model satisfies NLA there must exist $q \in \tilde{Q}$ such
that $q + \overline{q} \in \tilde{Q} \cup \{0\}$.  Suppose $q =
(m,\lambda)$ so that $q + \overline{q} = (2m,\lambda +
\overline{\lambda})$.  Let $q_0 = (0,\lambda)$.  We number the
four possibilities following the statement of the theorem.

(1)  $\lambda = \rho + i0$ with $\rho \neq 0$:  By (\ref{3.18}) $q =
(m,\rho)$ and $q + \overline{q} = (2m,2\rho)$ are in $Q$ which is
closed.  Hence, $q_0 = (0,\rho)$ and $q_0 + \overline{q}_0 =
(0,2\rho)$ are in $Q$.  $f^{q_0}(t) = \exp(\rho t)$ and $f^{q_0 + \overline{q}_0}(t) \ = \exp(2 \rho t)$
are in $Y(\R^n)$.

(2)  $\lambda = \rho + i \omega$ with $\rho \neq 0$ and $\omega
\neq 0$: Replacing $q$ by $\overline{q}$ if necessary we may
assume $\omega > 0$.  Again $q$ and $q + \overline{q} =
(2m,2\rho)$ are in $Q$, and again $q_0$ and $q_0 + \overline{q}_0
= (0,2\rho)$ are in $Q$.  Hence, Re$(f^{q_0}(t)) \ = \exp(\rho t)\cos \omega t$, Im$(f^{q_0}(t)) \ = \exp(\rho t) \sin \omega t$
and $f^{q_0 + \overline{q}_0}(t) \ = \exp(2 \rho t)$ are in $Y(\R^n)$.

(4)  $\lambda = 0 + i \omega$ with $\omega \neq 0$:  Again we can
assume $\omega > 0$ and as before $q$ and $q_0 = (0,i\omega)$ are
in $Q$.  Hence, Re$(f^{q_0}(t)) \ = \cos \omega t$ and
Im$(f^{q_0}(t)) \ = \sin \omega t$ are in $Y(\R^n)$.

(3)  $\lambda = 0 + i 0$:  Since $q = (m,0) \in \tilde{Q}$, $m \geq 1$. Furthermore, $q + \overline{q} = (2m,0) \in \tilde{Q}$. 
By (\ref{3.18}), $(m - 1,0)$ and $(2m - 1,0)$ are in $Q$.  Because $Q$ is
closed $(0,0)$ and $(1,0)$ are in $Q$.  Hence, $1 \ = t^0 \exp(0t)$
and $t \ = t^1 \exp(0t)$ are in $Y(\R^n)$.

To show that models which contain (1), (2) and (3) satisfy the NLA condition we
apply Lemma \ref{lem3.01} using the basis in each case given by Lemma \ref{lem3.04}.

(1)  With $q = (0,\rho),$ let $F = l^q$ so that $F(t) = \exp(\rho
t) - 1$. $ (F(t))^2 = (\exp(2 \rho t) - 1) - 2(\exp(\rho t) - 1)
= l^{q + \overline{q}}(t) - 2 l^q(t)$.  Hence, $F$, $(F)^2 \in
L(\R^n)$. Furthermore, $F: \R_+ \rightarrow \R_+$ is
injective.

(2)  With $q = (0,\rho + i \omega),$ let $F_1(t) = {\rm
Re}(l^q(t)) \ = \exp(\rho t)[\cos(\omega t)] - 1$ and $F_2(t) =
{\rm Im}(l^q(t)) \ = \exp(\rho t)[\sin(\omega t)]$.  $(F_1(t))^2 +
(F_2(t))^2 \ =  l^{q + \overline{q}}(t) - 2 F_1(t)$.  Hence, $F_1$,
$F_2$ and $(F_1)^2 + (F_2)^2$ are in $L(\R^n)$.  Since
$l^{q+\overline{q}}: \R_+ \rightarrow \R_+$ is
injective it follows that $(F^1,F^2): \R_+ \rightarrow \R^2$ is injective.

(3)  With  $(0,0)$ and $q = (1,0)$ in $Q$, $q$ and $q' = (2,0)$ are in $ \tilde{Q}$. Let $F(t) \ =
l^{q}(t) = t$ so that $(F(t))^2 \ = t^2 = l^{q'}(t)$.
  Hence, $F$ and $(F)^2$ are in $L(\R^n)$.
Clearly, $F: \R_+ \rightarrow \R_+$ is injective.

(4)  With $q = (0,i\omega)$, let $F_1(t) = {\rm Re}(l^q(t)) =
\cos(\omega t) - 1$ and $F_2(t) = {\rm Im}(l^q(t)) = \sin \omega
t$.  $(F_1(t))^2 + (F_2(t))^2 = -2 F_1(t)$.  Hence, $F_1$, $F_2$
and $(F_1)^2 + (F_2)^2$ are in $L(\R^n)$.  However,
$(F_1,F_2):\R_+ \rightarrow \R^2$ is not injective because $F_1$ and $F_2$ are
$2\pi/\omega$ periodic.

Thus, Lemma \ref{3.01} applies to cases (1), (2) and (3) but not to case
(4).

$\Box$ \vspace{.5cm}

{\bf Remark.}  If $\omega_1,\omega_2 > 0$ and the ratio
$\omega_1/\omega_2$ is irrational then with $q_\alpha =
(0,i\omega_\alpha)$ $(\alpha = 1,2)$ we define $F_1(t) + iF_2(t)
\ = l^{q_1}(t)$ and $F_3(t) + iF_4(t) \ = l^{q_2}(t)$.
$\sum^4_{k=1}\ (F_k(t))^2 \ = -2 F_1(t) - 2F_3(t)$.  Furthermore,
$(l^{q_1},l^{q_2}): \R_+ \rightarrow {\C}^2$ is
injective.  So Lemma \ref{lem3.01} shows that any linear yield curve model which
satisfies SPA and which contains the four parameter linear model
\begin{equation}\label{3.22}
Y(r_1,r_2,r_3,r_4) \ = \ r_1 \cos \omega_1t + r_2 \sin \omega_1t +
r_3 \cos \omega_2t + r_4 \sin \omega_2t
\end{equation}
satisfies NLA.  \vspace{.5cm}

%

{\bf Proof of Lemma \ref{lem3.05}:}  Define $Q^2 = \{q_1 + q_2: q_1,q_2 \in
\tilde{Q} \cup \{0\}\}$ and $Q^\# = \{q + \overline{q}: q \in
\tilde{Q}\}$.  Clearly, $Q^\#$ consists of real exponents only.
$Q^2$ and $\tilde{Q} \cup \{0\}$ are closed because $Q$ is
closed.  $Q^2$ contains both $Q^\#$ and $\tilde{Q} \cup \{0\}$.
The hypothesis of the lemma says exactly that $Q^\#$ and $\tilde{Q}
\cup \{0\}$ are disjoint, i.e.
\begin{equation}\label{3.23}
\tilde{Q} \cup \{0\} \subset Q^2\backslash Q^\#.
\end{equation}

Notice first that the assumptions imply that
\begin{equation}\label{3.23a}\begin{split}
(m, \ii \om) \not\in Q \quad \text{with} \ \  \om \not= 0,\\
(m,0) \not\in Q \quad \text{with} \ \ m \ge 1.
\end{split}\end{equation}
If $q = (m, \ii \om) \in Q$ then $q_0 = (0,\ii \om) \in \tilde Q$ and  $(0,0) = q_0 + \overline{q}_0 \in \tilde Q \cup \{ 0 \}$.
If $q = (m,0) \in Q$ with $m \ge 1$, then $(0,0),(1,0) \in Q$ and so $(1,0),(2,0) \in \tilde Q$. With $q_0 = (1,0) \in \tilde Q$,
$(2,0) = q_0 + \overline{q}_0 \in \tilde Q \cup \{ 0 \}$.
\vspace{.25cm}

By Lemma \ref{lem3.02}a, the set $\{f^q: q \in Q^2\}$ restricts to a linearly
independent set in the complex vector space ${\C}^T$ for $T$
any nonempty open set in $\R$.  By Corollary \ref{cor3.03} the exponent
set $Q$ for the $n$ parameter model has $n$ elements.  So $Q^2$
has at most $\frac{1}{2}(n+1)(n+2)$ elements.  By Lemma \ref{lem2.07} and the Remark thereafter,
the set $\{f^q(t) : t \in T \}$ spans the vector space ${\C}^{(Q^2)}$.
Thus, if we select for each $q \in Q^\#$ a positive
number $C_q$, then we can choose a finite function $z_t$ whose
support consists of at most $\frac{1}{2}(n+1)(n+2)$ times in $T$,
which we assume is contained in $(0,\infty)$, such that
\begin{equation}\label{3.24}
\begin{array}{ll}
\langle f^q(t)\rangle = 0 & q \in Q^2 \backslash Q^\#\\
\langle f^q(t)\rangle = C_q & q \in Q^\#.\end{array}
\end{equation}
Here we are using the notation of (\ref{3.03}) $\langle f(t) \rangle =
\sum_t z_t f(t)$.  It will be convenient to write $\langle f(t)
\rangle$ rather than $\langle f \rangle$ in these arguments.

A priori the numbers $z_t$ may be complex.  But $f^q(t)$ is a
real function for $q \in Q^\#$ and $\overline{f^q(t)} =
f^{\overline{q}}(t)$.  So $Q^2 - Q^\#$ closed under conjugation
implies that the conjugate vector $\overline{z}_t$ would yield
equations (\ref{3.24}) as well.  Hence, so would $\{{\rm Re}(z_t) =
\frac{1}{2}(z_t + \overline{z}_t)\}$.  Thus, we can choose $z_t$
to be a real vector.

Now let $\tilde{Q}_0 = \{q \in \tilde{Q}: \omega(q) = 0\}$ and
$\tilde{Q}_+ = \{q \in \tilde{Q}: \omega(q) > 0\}$.  From (\ref{3.23})
and (\ref{3.24}) we can take real and imaginary parts to get:
\begin{equation}\label{3.25} \begin{split}
\langle 1 \rangle = 0, \qquad
\langle f^q(t) \rangle = 0 \;\;\; q \in \tilde{Q}_0  \hspace{3cm}\\
\langle f^{\rho(q)}(t) \cos (\omega(q)t)\rangle =
\langle f^{\rho(q)}(t) \sin(\omega(q)t)\rangle = 0 \;\;\;\; q \in
\tilde{Q}_+.
\end{split}\end{equation}

Furthermore, we will now prove that

\[\langle (f^q(t))^2 \rangle = C_{2q} \;\;\;\; q \in
\tilde{Q}_0.\]
\[\langle[f^{\rho(q)}(t)\cos (\omega(q)t)]^2\rangle =
\langle[f^{\rho(q)}(t) \sin(\omega(q)t)]^2\rangle = \frac{1}{2}
C_{q+\overline{q},}\]
\begin{equation}\label{3.26}
\langle [f^{\rho(q)}(t) \cos
(\omega(q)t)][f^{\rho(q)}(t)\sin(\omega(q)t)]\rangle = 0 \;\;\; q
\in \tilde{Q}_+.
\end{equation}

For $q \in \tilde{Q}_0$, $2q = q + \overline{q} \in Q^\#$ and
$(f^q(t))^2 = f^{2q}(t)$.  So the first equation follows from
(\ref{3.24}) for $2q$.

For $q \in \tilde{Q}_+$, $q + \overline{q} = 2\rho(q) \in Q^\#$
and so $\sin^2 + \cos^2 = 1$ implies
\[\langle [f^{\rho(q)}(t) \cos (\omega(q)t)]^2\rangle + \langle
[f^{\rho(q)}\sin(\omega(q)t)]^2\rangle\]
\begin{equation}\label{3.27}
= \langle f^{q+\overline{q}}(t) \rangle = C_{q + \overline{q}}.
\end{equation}
On the other hand, since $q + q$ is not real it is in $Q^2
\backslash Q^\#$ and so $\langle f^{2q}(t)\rangle = 0$.  Taking real
and imaginary points we have
\begin{equation}\label{3.28}\begin{split}
0 = \langle f^{2\rho(q)}(t)\cos(2 \omega(q)t)\rangle = \hspace{5cm}\\
\langle[f^{\rho(q)}(t)\cos (\omega(q)t)]^2\rangle
- \langle [f^{\rho(q)}(t)\sin(\omega(q)t)]^2\rangle. \hspace{3cm}\\
\; \; 0 = \langle[f^{2\rho(q)}(t)\sin(2 \omega(q)(t)]\rangle
= 2\langle[f^{\rho(q)}(t)\cos
(\omega(q)t)][f^{\rho(q)}(t)\sin(\omega(q)t)]\rangle.
\end{split}\end{equation}
Together with (\ref{3.27}) these imply the rest of (\ref{3.26}).

By Lemma \ref{lem3.04} we can use as the basis $\{L^i: i = 1, \ldots, n\}$
for $L(\R^n)$ the set $\{l^q: q \in \tilde{Q}_0\} \cup \{
{\rm Re}(l^q),{\rm Im} (l^q): q \in \tilde{Q}_+\}$ listed in any
order. Recall from (\ref{3.19}) that $l^q(t) = f^q(t) - f^q(0)$. Thus, up to an
added constant, the basis   $\{L^i: i = 1, \ldots, n\}$ consists of the set
$\{f^q: q \in \tilde{Q}_0\} \cup \{
{\rm Re}(f^q),{\rm Im} (f^q): q \in \tilde{Q}_+\}$ listed in the same
order.
Thus (\ref{3.25}) implies $\langle 1 \rangle = 0$ and $\langle
L^i \rangle = 0$ for $ i = 1, \ldots, n$.  We show that for some
choice of positive constants $C_q$ the Hessian matrix $\langle
L^iL^j\rangle$ is positive definite.  Again in computing the Hessian we
may replace each $l^q$ by the corresponding $f^q$.

As an example, consider the case $Q = \{ (0,0) \}$, which is the flat yield curve example.
Here $\tilde Q = \{ (1,0) \}, Q^{\#} = \{ (2,0) \}$ and $Q^2$ is the disjoint union of
$Q^{\#}$ and $\tilde Q  \cup \{ 0 \}$. $L(\R)$ is one dimensional with basis $l^{(1.0)}(t) = t$.
Any positive choice for $C_{(2,0)}$ yields a positive definite $1 \times 1$ Hessian.

Returning to the general case, we see that the equations (\ref{3.26}) imply that the
diagonal entries are positive and that many of the off-diagonal
entries are zero.  However, they need not all be zero and this is why
we have to be a bit careful about the choice of $C_q$'s.

To illustrate this consider the case
\begin{equation}\label{3.29}
Q = \tilde{Q} = \{(0,-1),(0,-3),(0,-5)\}.
\end{equation}
Write $C_a$ for $C_{(0,-a)}$ $(a = 2, 6 , 10)$.  The Hessian
matrix in this case is
\begin{equation}
B = \left ( \begin{array}{ccc}
C_2 & 0 & C_6\\
0 & C_6 & 0\\
C_6 & 0 & C_{10}
\end{array} \right ).
\end{equation}

Now if we simply chose $C_a = 1$ for $a = 2,6,10$ then $B$ would
not be positive definite.  The remainder of the proof is
exemplified here by the idea that if $C_2 > 0$ and $C_6 > 0$ are
already chosen then $B$ is always positive definite provided that
$C_{10}$ is chosen large enough.

We complete the proof that $C_q$'s can always be chosen so that
the Hessian is positive definite by using induction on the number of
points in $Q$.

Begin by choosing a set of extreme points $Q^*$ of $Q$ as
follows.  Among the nonnegative numbers $\{|{\rm Re}(\lambda)|:
(m,\lambda) \in Q\}$ choose the largest and denote it by $a^*$.
If $a^* = 0$, then (\ref{3.23a}) implies that $Q = \{ (0,0) \}$, the
flat yield curve case we considered above, and so we may assume that $a^* > 0$.

At least one of the sets $\{(m,\lambda) \in Q: {\rm Re}(\lambda)
= a^*\}$ and $\{(m,\lambda) \in Q: {\rm Re}(\lambda) = -a^*\}$ is
nonempty.  We will suppose that the first one is nonempty and let
$m^* = \max\{m: (m,\lambda) \in \tilde{Q}$ and ${\rm Re}(\lambda)
= a^*\}$.  Now let $Q^* = \{(m^*,\lambda) \in \tilde{Q}:
{\rm Re}(\lambda) = a^*\}$.  Since $a^* > 0$, $Q^{*} \subset Q$. Notice that $q^* \in Q^*$
implies $q^* + \overline{q}^* = (2m^*,2a^*)$ which we will denote
by $\rho^*$.

Define $Q^\prime = Q \setminus Q^*$ and notice that
$Q^\prime$ is a smaller closed subset still satisfying the
hypothesis of the Theorem, i.e. $\tilde{Q}^\prime \cup \{0\}
\subset Q^{\prime 2} \setminus Q^{\prime \#}$.  Furthermore, $\rho^*
\not\in Q^{\prime \#}$.  In fact, we have
\begin{equation}\label{3.30}
\rho^* \not\in \{{\rm Re}(q_1 + q_2): q_1 \in \tilde{Q}^\prime
\;\;{\rm and} \;\; q_2 \in \tilde{Q} \}.
\end{equation}
For suppose $q_1 = (m_1,\lambda_1)$ and $q_2 = (m_2,\lambda_2)$.
By construction either ${\rm Re}(\lambda_1) < a^*$ or $m_1 < m^*$. Furthermore,
${\rm Re}(\lambda_2) \leq a^*$ and if ${\rm Re}(\lambda_2) = a^*$, then $m_2 \leq m^*$.
  Hence, ${\rm Re}(\lambda_1 + \lambda_2) \le 2a^*$ and if ${\rm Re}(\lambda_1 + \lambda_2) = 2a^*$,
  then $m_1 + m_2 < 2m^*$.

Now arrange the list $\{L^i: i = 1, \ldots, n\}$ so that the
functions $l^q$ with $q \in Q^*$ occur at the end.  The
Hessian matrix $B$ has the bloc form.
\begin{equation}\label{3.31}
B = \left (\begin{array}{cc}
B^\prime & U\\
U^T & B^*
\end{array} \right )
\end{equation}
where $B^\prime$ is the matrix associated with $Q^\prime$.   By
inductive hypothesis we can choose $C_q > 0$ for $q \in
Q^{\prime \#}$ so that $B^\prime$ is positive definite.

Now some of the $C_q$'s or $\frac{1}{2}C_q$'s may also appear in
$U$ and $U^T$ for $q \in \tilde{Q}^\prime$.  This is what
happened in the above example.  However, (\ref{3.30}) implies that $q_1 +
q_2$ in $Q^2$ can equal $\rho^*$ only for $q_1,q_2$ in
$\tilde{Q}^*$.  Thus, if the $C_q$'s have been chosen for $q
\in Q^{\prime 2}$ so that $B^\prime$ is positive definite, then
$U$ is also determined.

Furthermore, $B^*$ is a positive diagonal matrix.  By (\ref{3.26})
every diagonal entry is either $C_{\rho^*}$ or $\frac{1}{2}
C_{\rho^*}$.  If $Q^*$ consists of the real exponent
$(m^*,a^*)$ alone then $B^*$ is the $1 \times 1$ matrix
$(C_{\rho^*})$.  Otherwise the off-diagonal entries are all zero
because if $q_1,q_2 \in \tilde{Q}^*$ is not a conjugate pair then
$q_1 + q_2$ is not real and so is not in $Q^\#$.  By (\ref{3.24}),
$\langle f^{q_1+q_2}\rangle$ and its conjugate are zero.

In detail, suppose that $q_1 = (m^*,a^* + \ii \omega_1), q_2 = (m^*,a^* + \ii \omega_2)$ so that
$q_1 + q_2 =  (2m^*,2a^* + \ii (\omega_1 + \omega_2)), q_1 + \overline{q}_2 =  (2m^*,2a^* + \ii (\omega_1 - \omega_2))$.
It follows that $\langle f^{2 \rho^*}(t)\cos(\omega_1 + \omega_2)t \rangle = 0 = \langle f^{2 \rho^*}(t)\cos(\omega_1 - \omega_2)t \rangle$
and $\langle f^{2 \rho^*}(t)\sin(\omega_1 + \omega_2)t \rangle = 0 = \langle f^{2 \rho^*}(t)\sin(\omega_1 - \omega_2)t \rangle$.
Now apply the sum formulae for sine and cosine.

It follows that $B^* =
C_{\rho^*} \Delta$ where $\Delta$ is a diagonal matrix with diagonal
entries all $1$ or $\frac{1}{2}$.

In the special case where $Q = Q^*$, $B = B^*$ and any positive
choice of $C_{\rho^*}$ will do.  This includes the initial step
for the induction.

In general $B$ is positive definite provided $C_{\rho^*}$ is now
chosen large enough.  This comes from the following lemma whose
proof completes the proof of Lemma \ref{lem3.05}.

\vspace{2ex}

\begin{lem}\label{lem3.06} Let $B^\prime$ be a $k \times k$
positive definite symmetric matrix, $D$ be a $n - k \times n - k$
positive definite symmetric matrix and $U$ be a $k \times n - k$
matrix.  If $C > 0$ is chosen sufficiently large then
\[ B = \left ( \begin{array} {cc}
B^\prime & U\\
U^T & CD
\end{array} \right )\]
is positive definite.\end{lem}

\proof Write $z$ in ${\R}^n$ as $z = (x,y)$ where $x
\in {\R}^k$ and $y \in {\R}^{n-k}$.  Then
\[z^TB z = x^TB^\prime x + 2x^TUy + Cy^TDy.\]

We want to choose $C > 0$ large enough that $z \neq 0$ implies
$z^TBz > 0$.

Because $B^\prime$ and $D$ are positive definite we can find
$\epsilon > 0$ and $M > 0$ such that
\newpage
\[x^TB^\prime x \geq \epsilon\parallel x\parallel^2\]
\[y^TDy \geq \epsilon\parallel y\parallel ^2\]
\[|2x^TUy| \leq M\parallel x  \parallel \parallel y \parallel.\]

Now if $\epsilon M^{-1} \parallel x \parallel > \parallel y
\parallel$ then

\[z^TB z \geq x^TB^\prime x - |2x^TUy|\]
\[\;\;\;\;\;> \epsilon \parallel x \parallel^2 - M \parallel x
\parallel (\epsilon M^{-1}\parallel x \parallel) = 0.\]

Similarly, $C\epsilon M^{-1} \parallel y \parallel > \parallel x
\parallel$ implies
\[z^TB z > Cy^TDy - |2x^TUy|\]
\[\;\;\;\;\; > C \epsilon \parallel y \parallel^2 - M(C \epsilon
M^{-1}\parallel y \parallel) \parallel y \parallel = 0.\]

If we choose $C > (M/\epsilon)^2$ then for any $z \neq 0$ at
least one of the two inequalities will hold, i.e. if $\parallel y
\parallel\; \geq \epsilon M^{-1} \parallel x \parallel$ and
$\parallel y \parallel \; > 0$ then $C \epsilon M^{-1}\parallel y
\parallel \; > M \epsilon^{-1}\parallel y \parallel \ \geq \ \parallel x \parallel$.

$\Box$  \vspace{.5cm}

{\bf Proof of Lemma \ref{lem3.01}:}  For $j = 1, \ldots, p$ there exist real
constants $u^j_i $ ($i = 1, \ldots , n)$ such that $F^j = \sum_i
u_i^j L^i$ and $u_i^{p+1}$ ($i = 1, \ldots, n$) such that
$\sum_j(F^j)^2 = \sum_i u_i^{p+1}L^i$.  For $s \in {\R}^p$ and
$S \in {\R}$ define $r(s,S)$ by $r(s,S)_i = -\sum_j s_ju^j_i -
S u^{p+1}_i$ ($i = 1, \ldots, n)$.

Suppose we are given a finite function $z_t$ on $\R_+$ such
that $z_0 = 0$ and $\langle 1 \rangle = \langle L^i \rangle = 0$
for $i = 1, \ldots, n$.  Clearly,
\[\langle 1 \rangle = 0\]
\[\langle F^j \rangle = 0 \;\;\; (j = 1, \ldots, p)\]
\begin{equation}\label{3.32}
\langle \sum^p_{j=1}(F^j)^2 \rangle = 0.
\end{equation}
Now given $r^*$ we define the basket of futures $x_t$ by (\ref{3.04}), i.e. $x_t = z_t/P_t(r^*)$.
Then from (\ref{3.06}) we have for $r = r^* + r(s,S)$ that $V(r)$ is
given by
\[V(s,S) \equiv \langle \exp [\sum^p_{j=1} s_jF^j + S
\sum^p_{j=1}(F^j)^2]\rangle \]
\begin{equation}\label{3.33}
= \langle \exp [s^TF + SF^TF]\rangle.
\end{equation}
For the last expression we are regarding $s$ and each $F_t$ as
column vectors in ${\R}^p$.  By (\ref{3.05}) and (\ref{3.07}), $x_t$ is a
self-financing basket of futures with respect to the original
price curve $P_t(r^*)$ and $r = r^*$ is a critical point for the
value function $V(r)$.  In particular, at the origin $(s,S) =
(0,0)$ $V(0,0)$ equals zero.  For $s$ a nonzero vector in ${\R}^p$
define
\begin{equation}\label{3.34}
\epsilon = \parallel s \parallel = (s^Ts)^{1/2} \;\;{\rm and }
\;\; \theta = s/\parallel s \parallel.
\end{equation}
Thus, $\theta$ is a variable on the unit sphere in ${\R}^p$.

For $\epsilon \neq 0$ we define $M = S/\epsilon^2$.  So we have
the change of variables:
\begin{equation}\label{3.35}
s = \epsilon \theta \;\;\;{\rm and} \;\;\; S = M \epsilon^2.
\end{equation}

We prove there exist $(\theta_+,M_+)$ and $(\theta_-,M_-)$ such
that for $\epsilon > 0$ small enough
\[V(\epsilon \theta_+,\epsilon^2M_+) > 0\]
\begin{equation}\label{3.36}
V(\epsilon \theta_-,\epsilon^2 M_-) < 0.
\end{equation}
It follows that the origin is not a local minimum (or local maximum)
point for $V(s,S)$ and so $r^*$ is not for $V(r)$.

Notice first that $V(s,S)$ does not vanish identically on any
neighborhood $O$ of the origin.  For if $\sum z_t \exp [s^TF_t +
SF^T_tF_t] = 0$ as a function on $O$ then by the Remark after Lemma
\ref{lem3.02} it would have to be true that for some $t_1 \neq t_2$ both in
the support of $\{z_t\}$.  $F^j_{t_1} = F^j_{t_2}$ ($j = 1,
\ldots, p$).  This contradicts the assumption of condition (1) in
the statement of Lemma \ref{lem3.01}. That is, by that assumption $(F^1_{t_1},\dots,F^p_{t_1}) \not= (F^1_{t_2},\dots,F^p_{t_2})$ for
for $t_1 \not= t_2$ in the support of $z$ and so the functions $\exp [s^TF_t + SF^T_tF_t]$
of $(s,S) \in O$ are linearly independent and we would not have $\sum z_t \exp [s^TF_t + SF^T_tF_t] = 0$ on $O$.

Now we use the series expansion for the exponential and the
binomial theorem.
\[V(s,S) = \sum^\infty_{n=1} \frac{1}{n!} \langle(s^TF +
SF^TF)^n\rangle\]
\[= \sum^\infty_{n=1} \sum^n_{k=0} \frac{1}{k!(n-k)!}\langle
(s^TF)^{n-k}S^k(F^TF)^k\rangle\]
\begin{equation}\label{3.37}
= \sum^\infty_{n=1} \sum^n_{k=0} \frac{M^k}{k!(n-k)!}
\epsilon^{n+k}\langle(\theta^TF)^{n-k}(F^TF)^k\rangle
\end{equation}
Observe that $\langle 1 \rangle = 0$ implies we can begin the
series with the $n = 1$ term.  Having used the polar coordinate
change (\ref{3.36}) we similarly write, with $\parallel F \parallel =
(F^TF)^{1/2}$,
\begin{equation}\label{3.38}
F = \parallel F \parallel \omega, \;\;\;{\rm i.e.} \;\;\; \omega
= F/\parallel F \parallel \;\; (\parallel F \parallel\; \neq 0),
\end{equation}
and change variables in the sum letting $N = n + k$:
\[V(s,S) = \sum^\infty_{n=1} \sum^n_{k=0} \frac{M^k}{k!(n-k)!}
\epsilon^{n+k}\langle \parallel F \parallel^{n+k}(\theta^T
\omega)^{n-k}\rangle\]
\begin{equation}\label{3.39}
= \sum^\infty_{N=1} \epsilon^N \sum^{[N/2]}_{k=0} \frac{M^k}{k!(N-2k)!}
\langle \parallel F \parallel^N(\theta^T \omega)^{N-
2k}\rangle,
\end{equation}
where $[N/2]$ is the greatest integer less than or equal to
$N/2$.

 Because this series does not vanish identically there is a
smallest positive value $N = N^*$ such that
\begin{equation}\label{3.40}
V_{N^*}(\theta,M) \equiv \sum^{[N^*/2]}_{k=0}
\frac{M^k}{k!(N^*-2k)!} \langle \parallel F
\parallel^{N^*}(\theta^T \omega)^{N^*-2k}\rangle
\end{equation}
is not identically zero.  Furthermore, absolute convergence of
the series implies for any $(\theta,M)$ such that
$V_{N^*}(\theta,M) \neq 0$ the sign of $V(\epsilon
\theta,\epsilon^2M)$ is the same as the sign of
$V_{N^*}(\theta,M)$ provided $\epsilon > 0$ is sufficiently
small.  Thus, it suffices to find $(\theta_+,M_+)$ and
$(\theta_-,M_-)$ such that
\[V_{N^*}(\theta_+,M_+) > 0\]
\begin{equation}\label{3.41}
V_{N^*}(\theta_-,M_-) < 0.
\end{equation}
Then (\ref{3.36}) follows for $\epsilon > 0$ sufficiently small.

We proceed by averaging the $\theta$ variable over the unit
sphere in ${\R}^p$.  There are two cases:
\vspace{.5cm}

{\bf Case i:}  $\int_{S^{p-1}}V_{N^*}(\theta,M)d\theta $ vanishes
for all $M$.
\vspace{.5cm}

In this case, there is some $M_0$ such that $V_{N^*}(\theta,M_0)$
is not identically zero as a function of $\theta$ by choice of
$N^*$.  Since the average in $\theta$ is $0$ at $M_0$, there must
exist values $\theta_+$ and $\theta_-$ satisfying (\ref{3.41}) with
$M_+ = M_0 = M_-$.
\vspace{.5cm}

{\bf Case ii:}  $\int_{S^{p-1}} V_{N^*}(\theta,M)d\theta$
does not vanish identically.
\vspace{.5cm}

In this case we compute the average.  The results depend on
whether $p = 1$ or $p > 1$.

When $p = 1$, $\theta$ and $\omega$ are each $\pm 1$.  For each
$\omega$, the function $\theta^T\omega$ takes the two values $\pm
1$, each weighted $1/2$.  So the average of $(\theta^T
\omega)^{N^*-2k}$ is $0$ when $N^*$ is odd and is $1$ when $N^*$
is even.

Thus, when $p = 1$, we have
\begin{equation}\label{3.42}
\int\; V_{N^*}(\theta ,M)d\theta = \left \{
\begin{array}{ll} 0 & N^* \;\;\;{\rm odd}\\
P_R(M) \cdot \langle \parallel F \parallel^{N^*}\rangle & N^* = 2R\\
\end{array} \right.
\end{equation}
where
\begin{equation}\label{3.43}
P_R(M) \ \equiv \ \sum^R_{k=0} \  \frac{M^k}{k!(2R - 2k)!}.
\end{equation}

When $p > 1$, we let $\alpha$ denote the angle between $\omega$,
fixed on $S^{p-1}$, and $\theta$.  Thus,
\begin{equation}\label{3.44}
\theta = (\cos \alpha) \omega + (\sin \alpha) \psi
\end{equation}
where $\psi$ varies on the $p-2$ sphere equatorial with respect
to the pole at $\omega$.  We normalize so that
\begin{equation}\label{3.45}
1 \ = \int_{S^{p-1}}1 d\theta \ = \int_{S^{p-2}}\int^\pi_0 d\alpha
d\psi \ = \pi \int_{S^{p-2}} d\psi.
\end{equation}
It follows that
\[\int_{S^{p-1}}(\theta^T \omega)^{N^*-2k}d \theta \ =
\int_{S^{p-2}}\int^\pi_0(\cos \alpha)^{N^*-2k}d\alpha
d\psi\]
\begin{equation}\label{3.46}
 \ = \frac{1}{\pi}\int^\pi_0 (\cos
\alpha)^{N^*-2k} d\alpha \ = \ \frac{1}{\pi} \ \left \{ \begin{array} {cl}
0 & N^* \;\;{\rm odd}\\
\frac{(N^*-2k-1)!!}{(N^*-2k)!!} & N^*\;\;{\rm even}
\end{array} \right. .
\end{equation}
The formula for $\int^\pi_0(\cos \alpha)^R
d\alpha$ is proved by using induction on $R$ and the formula
\begin{equation}\label{3.47}
\frac{d^2}{d \alpha^2} (\cos \alpha)^R = R(R - 1)(\cos
\alpha)^{R-2} - R^2 (\cos \alpha)^R.
\end{equation}
The {\em semi-factorial} $R!!$ for a positive integer $R$ is $R
\cdot (R - 2)\ldots$ descending to $1$ or $2$ .  In particular,
$(2R - 2k - 1)!!/(2R - 2k)! \ = 1/(2R - 2k)!!$.

Thus, when $p > 1$, we have
\begin{equation}\label{3.48}
\int V_{N^*}(\theta,M)d \theta \ = \left \{\begin{array} {ll} 0 &
N^* \;\;{\rm odd}\\
\tilde{P}_R(M) \langle \parallel F \parallel^{N^*}\rangle / \pi & N^* \ =
2R
\end{array} \right. .
\end{equation}
where
\begin{equation}\label{3.49}
\tilde{P}_R(M) \ \equiv \ \sum^R_{k=0} \ \frac{M^k}{k![(2R-2k)!!]^2}.
\end{equation}

Thus, for any positive integer $p$, if $N^*$ is odd or $\langle
\parallel F \parallel^{N^*} \rangle = 0$, then we return to case $i$.

When $N^* = 2R$ and $\langle \parallel F \parallel^{N^*}\rangle
\neq 0$ we complete the proof by using Lemma 7 below.  
When $p =
1$ it says that we can choose $M_+$ and $M_-$ so that
\begin{equation}\label{3.50}
P_R(M_+)\langle \parallel F \parallel^{N^*}\rangle > 0 \; {\rm
and} \; P_R(M_-) \langle \parallel F \parallel^{N^*}\rangle < 0.
\end{equation}
When the average of a function of $\theta$ is nonzero there will
certainly exist values having the same sign as the average.  So
by (\ref{3.42}) we can choose $\theta_+$ and $\theta_-$ to satisfy
(\ref{3.41}).  When $p > 1$ we use the same argument, replacing $P_R$
by $\tilde{P}_R$ in (\ref{3.50}) and replacing (\ref{3.42}) by (\ref{3.48}).  So the proof of
the following lemma completes the proof of Lemma \ref{lem3.01}.

\begin{lem}\label{lem3.07}For every positive integer $R$
define the polynomials of degree $R$:
\[P_R(M) = \sum^R_{k=0} \frac{M^k}{k!(2R-2k)!}\]
\[\tilde{P}_R(M) = \sum^R_{k=0} \frac{M^k}{k![(2R - 2k)!!]^2}.\]
If $M \geq 0$ then $P_R(M)$, $\tilde{P}_R(M) > 0$ and for every
$R$ there exist negative numbers $M$, $\tilde{M}$ such that
$P_R(M) < 0$ and $\tilde{P}_R(\tilde{M}) < 0$.\end{lem}

\proof Because the coefficients are positive it is clear
that $P_R(M) \geq P_R(0) > 0$ when $M > 0$.  Furthermore if $R$
is odd then ${\rm Lim}_{M\rightarrow - \infty} P_R(M) = -\infty$
and so $P_R(M) < 0$ for $M < 0$ with $|M|$ sufficiently large.
We complete the proof for $P_R$ by finding a range of negative
values $M$ where $P_R(M)$ has the same sign as the $R-1$ term of
the sum.  When $R$ is even this term is negative.

Let $T_k$ denote the $k^{th}$ term and $S_k$ the $k^{th}$ partial
sum: \begin{equation}\label{3.51}
T_k = \frac{M^k}{k!(2R-2k)!}
\end{equation}

\begin{equation}\label{3.52}
S_k = \sum^k_{j=0} T_j = S_{k-1} + T_k.
\end{equation}

Define the absolute ratio $(k = 1,\ldots, R)$:
\begin{equation}\label{3.53}
\frac{|T_{k-1}|}{|T_k|} = \frac{k}{(2R-2k+2)(2R-2k+1) |M|} =
\mu_k.
\end{equation}
Notice that $\mu_k$ increases monotonically with $k$ provided
$|M|$ is fixed.  If $|M| > (R - 1)/6$, then $\mu_{R-1} < 1/2$
(direct computation) implies $\mu_k < 1/2$ for all $k = 1,
\ldots, R - 1$.

By induction on $k$ we then prove that with $M$ negative and $|M|
> (R-1)/6$:
\begin{equation}\label{3.54}
\frac{1}{2} < \frac{S_k}{T_k} \leq 1 \;\;\;\; k = 0, \ldots, R -1.
\end{equation}
 Because $S_0 = T_0 \neq 0$ the initial step
is clear.  Inductively,
\begin{equation}\label{3.55}
\frac{S_k}{T_k} = \frac{S_{k-1} + T_k}{T_k} = 1 + \frac{S_{k-
1}}{T_k}.
\end{equation}
Because $M$ is negative the terms of the series alternate in sign
and so $T_{k-1}/T_k = -\mu_k$.  Thus,
\begin{equation}\label{3.56}
\frac{S_k}{T_k} = 1 - \mu_k \frac{S_{k-1}}{T_{k-1}}.
\end{equation}
Now $\mu_k > 0$ and $S_{k-1}/T_{k-1} > 0$ (in fact $> 1/2$) imply
$S_k/T_k < 1$.  Also, $\mu_k < 1/2$ and $S_{k-1}/T_{k-1} \leq 1$
imply $S_k/T_k > 1/2$.

Now choose $M$ negative and so that $(R-1)/6 < |M| < R/4$.  $S_R$
satisfies
\begin{equation}\label{3.57}
S_R = S_{R-1} + T_R = S_{R-1} - \frac{1}{\mu_R} T_{R-1}
\end{equation}
and so
\begin{equation}\label{3.58}
\frac{S_R}{T_{R-1}} = \frac{S_{R-1}}{T_{R-1}} - \frac{1}{\mu_R}.
\end{equation}

By (\ref{3.55}) $S_{R-1}/T_{R-1} > 1/2$.  On the other hand $|M| < R/4$
implies by direct computation again that $\mu_R > 2$ and so
$1/\mu_R < 1/2$.  So with $M$ negative in the range $(R-1)/6 < |M| < R/4$
\begin{equation}\label{3.59}
\frac{S_R}{T_{R-1}} > 0.
\end{equation}

That is, $P_R(M) = S_R$ has the same sign as $T_{R-1}$.  This sign is
negative when $R$ is even.

For $\tilde{P}_R(M)$ the proof follows exactly the same pattern.
Equation (\ref{3.53}) is replaced by
\begin{equation}\label{3.60}
\frac{|T_{k-1}|}{|T_k|} = \frac{k}{(2R - 2k + 2)^2 |M|} \equiv
\mu_k.
\end{equation}
The required range for $\tilde{M}$ negative such that
$\tilde{P}_R(\tilde{M})$ has the same sign as its $R-1$ term as above is
$(R-1)/8 < |\tilde{M}| < R/8$.

$\Box$ \vspace{1cm}

\newpage

\bibliographystyle{amsplain}

\begin{thebibliography}{10}

\bibitem{FW} L. Fisher and R. Weil, ``Coping with the Risk of Interest
Rate Fluctuations: Returns to Bondholders from Naive and Optima
Strategies'', J. of Business  (1971), \textbf{44}:408-431.
\vspace{.5cm}

\bibitem{H} J. Hicks, ``Value and Capital,'' Oxford U. Press, Oxford,
1939.
\vspace{.5cm}

\bibitem{ISW} J. Ingersoll, Jr., J. Skelton, and R. Weil, ``Duration Forty
Years Later'' J. of Financial and Quantitative Analysis, (1978),
627-650.
\vspace{.5cm}

\bibitem{L}  S. Lang, {\em Differential and Riemannian Manifolds}, (1995)
Springer-Verlag, Berlin and New York.
\vspace{.5cm}

\bibitem{M} F. Macauley, ``Some Theoretical Problems Suggested by the
Movement of Interest Rates, Bond Yields, and Stock Prices in the
United States Since 1856'', Columbia U. Press, New York, 1938.
\vspace{.5cm}

\bibitem{W} R. Weil, ``Macauley's Duration: An appreciation'', J. of
Business (1973), \textbf{46}:589-592.
\end{thebibliography}

\end{document}